\documentclass[draft]{amsart}

\usepackage{amsthm}
\usepackage{amsmath}
\usepackage{amssymb}
\usepackage{amsfonts}
\usepackage{latexsym}
\usepackage[all]{xy} \SelectTips{eu}{}
\usepackage[hidelinks=true]{hyperref}
\usepackage{xcolor}

\newcommand{\numberseries}{\bfseries}   
\newlength{\thmtopspace}                
\newlength{\thmbotspace}                
\newlength{\thmheadspace}               
\newlength{\thmindent}                  
\newtheoremstyle{bfupright head,slanted body}
                {\thmtopspace}{\thmbotspace}
                {\slshape}{\thmindent}{\bfseries}{.}{\thmheadspace}
                {{\numberseries \thmnumber{#2\;}}\thmnote{#3}}

\newtheoremstyle{bfupright head,upright body}
                {\thmtopspace}{\thmbotspace}
                {\upshape}{\thmindent}{\bfseries}{.}{\thmheadspace}
                {{\numberseries \thmnumber{#2\;}}\thmnote{#3}}

\newtheoremstyle{fixed bf head,slanted body}
                {\thmtopspace}{\thmbotspace}{\slshape}
                {\thmindent}{\bfseries}{.}{\thmheadspace}
                {{\numberseries \thmnumber{#2\;}}\thmname{#1}\thmnote{ (#3)}}

\newtheoremstyle{fixed bf head,upright body}
                {\thmtopspace}{\thmbotspace}{\upshape}
                {\thmindent}{\bfseries}{.}{\thmheadspace}
                {{\numberseries \thmnumber{#2\;}}\thmname{#1}\thmnote{ (#3)}}

\newtheoremstyle{numbered paragraph}
                {\thmtopspace}{\thmbotspace}{\upshape}
                {\thmindent}{\upshape}{}{\thmheadspace}
                {{\numberseries \thmnumber{#2.}}}

\theoremstyle{bfupright head,slanted body}
\newtheorem{res}{}[section]             \newtheorem*{res*}{}

\theoremstyle{bfupright head,upright body}
               \newtheorem*{bfhpg*}{}

\theoremstyle{fixed bf head,slanted body}
\newtheorem{thm}[res]{Theorem}          \newtheorem*{thm*}{Theorem}
\newtheorem{prp}[res]{Proposition}      \newtheorem*{prp*}{Proposition}
\newtheorem{cor}[res]{Corollary}        \newtheorem*{cor*}{Corollary}
\newtheorem{lem}[res]{Lemma}            \newtheorem*{lem*}{Lemma}

\theoremstyle{fixed bf head,upright body}
\newtheorem{dfn}[res]{Definition}       \newtheorem*{dfn*}{Definition}
\newtheorem{rmk}[res]{Remark}           \newtheorem*{rmk*}{Remark}
\newtheorem{exa}[res]{Example}           \newtheorem*{exa*}{Example}
           \newtheorem*{que*}{Question}
\newtheorem{fct}[res]{Fact}           \newtheorem*{fact*}{Fact}

\theoremstyle{numbered paragraph}

\newlength{\thmlistleft}        
\newlength{\thmlistright}       
\newlength{\thmlistpartopsep}   
\newlength{\thmlisttopsep}      
\newlength{\thmlistparsep}      
\newlength{\thmlistitemsep}     

\newcounter{eqc}
\newenvironment{eqc}{\begin{list}{\upshape (\textit{\roman{eqc}})}%
    {\usecounter{eqc}%
      \setlength{\leftmargin}{\thmlistleft}%
      \setlength{\labelwidth}{\thmlistleft}%
      \setlength{\rightmargin}{\thmlistright}%
      \setlength{\partopsep}{\thmlistpartopsep}%
      \setlength{\topsep}{\thmlisttopsep}%
      \setlength{\parsep}{\thmlistparsep}%
      \setlength{\itemsep}{\thmlistitemsep}}}%
  {\end{list}}%

\newcounter{prt}
  {\end{list}}%

\newcounter{rqm}
  {\end{list}}%

\newenvironment{prf*}[1][Proof]{%
  \begin{proof}[\bf #1]
    \setcounter{equation}{0}
    }
  {\end{proof}
}

\newcommand{\pgref}[1]{\ref{#1}}
\newcommand{\thmref}[2][Theorem~]{#1\pgref{thm:#2}}
\newcommand{\corref}[2][Corollary~]{#1\pgref{cor:#2}}
\newcommand{\prpref}[2][Proposition~]{#1\pgref{prp:#2}}
\newcommand{\lemref}[2][Lemma~]{#1\pgref{lem:#2}}
\newcommand{\dfnref}[2][Definition~]{#1\pgref{dfn:#2}}
\newcommand{\rmkref}[2][Remark~]{#1\pgref{rmk:#2}}
\newcommand{\fctref}[2][Fact~]{#1\pgref{fct:#2}}
\newcommand{\secref}[2][Section~]{#1\pgref{sec:#2}}
\renewcommand{\eqref}[1]{(\pgref{eq:#1})}

\newcommand{\thmcite}[2][?]{\cite[Thm.~#1]{#2}}
\newcommand{\prpcite}[2][?]{\cite[Prop.~#1]{#2}}
\newcommand{\corcite}[2][?]{\cite[Cor.~#1]{#2}}
\newcommand{\lemcite}[2][?]{\cite[Lem.~#1]{#2}}
\newcommand{\seccite}[2][?]{\cite[Sec.~#1]{#2}}
\newcommand{\dfncite}[2][?]{\cite[Def.~#1]{#2}}
\newcommand{\rmkcite}[2][?]{\cite[Rmk.~#1]{#2}}
\newcommand{\exacite}[2][?]{\cite[Exa.~#1]{#2}}
\newcommand{\eqclbl}[1]{{\upshape(\textit{#1})}}
\newcommand{\proofofimp}[3][:]{\mbox{\eqclbl{#2}$\!\implies\!$\eqclbl{#3}#1}}

\newcommand{\Aop}{{A^\circ}}

\numberwithin{equation}{res}
\def\urltilda{\kern -.15em\lower .7ex\hbox{\~{}}\kern .04em}

\newcommand{\setof}[3][\mspace{1mu}]{\{#1#2 \mid #3#1\}}

\newcommand{\ZZ}{\mathbb{Z}}
\newcommand{\QQ}{\mathbb{Q}}
\newcommand{\NN}{\mathbb{N}}
\newcommand{\xra}[2][]{\xrightarrow[#1]{\;#2\;}}
\newcommand{\qra}{\xra{\qis}}
\newcommand{\qis}{\simeq}
\newcommand{\Hom}[1][A]{\operatorname{Hom}_{#1}}
\newcommand{\Ext}{\operatorname{Ext}}

\newcommand{\id}[1][A]{\operatorname{id}_#1}
\newcommand{\pd}[1][A]{\operatorname{pd}_#1}
\newcommand{\fd}[1][A]{\operatorname{fd}_{#1}}
\newcommand{\fpid}[1][A]{\operatorname{fp-id}_{#1}}
\newcommand{\Gid}[1][A]{\operatorname{Gid}_#1}
\newcommand{\Gpd}[1][A]{\operatorname{Gpd}_#1}
\newcommand{\Coker}{\operatorname{Coker}}
\newcommand{\Ggldim}[1][A]{\operatorname{Ggldim}(#1)}
\newcommand{\Gwgldim}[1][A]{\operatorname{Gwgldim}(#1)}
\newcommand{\Gfcgldim}[1][A]{\operatorname{Gfc-gldim}(#1)}
\newcommand{\splf}[1][A]{\operatorname{splf}(#1)}
\newcommand{\silf}[1][A]{\operatorname{silf}(#1)}
\newcommand{\silp}[1][A]{\operatorname{silp}(#1)}
\newcommand{\spli}[1][A]{\operatorname{spli}(#1)}
\newcommand{\sfli}[1][A]{\operatorname{sfli}(#1)}
\newcommand{\FFD}[1][A]{\operatorname{FFD}(#1)}
\newcommand{\FID}[1][A]{\operatorname{FID}(#1)}
\newcommand{\FPD}[1][A]{\operatorname{FPD}(#1)}
\newcommand{\Gfcd}[1][A]{\operatorname{Gfcd}_#1}
\newcommand{\Gfd}[1][A]{\operatorname{Gfd}_#1}
\newcommand{\Inj}[1][A]{\mathsf{Inj}(#1)}
\newcommand{\Prj}[1][A]{\mathsf{Prj}(#1)}
\newcommand{\Mod}[1][A]{\mathsf{Mod}(#1)}
\newcommand{\GI}{\mathsf{GInj}(A)}
\newcommand{\GP}{\mathsf{GPrj}(A)}
\newcommand{\Flat}{\mathsf{Flat}(A)}
\newcommand{\Cot}{\mathsf{Cot}(A)}
\newcommand{\SC}{\mathsf{SCot}(A)}
\newcommand{\FpInj}[1][A]{\mathsf{FpInj}(#1)}
\newcommand{\sfpilp}[1][A]{\operatorname{fp-silp}(#1)}
\newcommand{\FFPID}[1][A]{\operatorname{fp-FID}(#1)}
\newcommand{\lra}{\longrightarrow}
\newcommand{\Filt}{\operatorname{Filt}}
\newcommand{\add}{\operatorname{add}}

\newcommand{\wrightperiodic}{weakly right periodic}
\newcommand{\wleftperiodic}{weakly left periodic}
\newcommand{\rightperiodic}{right periodic}
\newcommand{\leftperiodic}{left periodic}

\newcommand{\rightbothperiodic}{(weakly) right periodic}
\newcommand{\leftbothperiodic}{(weakly) left periodic}

\newcommand{\sfL}{\mathsf{L}}
\newcommand{\sfC}{\mathsf{C}}
\newcommand{\sfS}{\mathsf{S}}
\newcommand{\sfU}{\mathsf{U}}
\newcommand{\sfV}{\mathsf{V}}
\newcommand{\sfW}{\mathsf{W}}
\newcommand{\sfUV}{(\sfU,\sfV)}
\newcommand{\capUV}{\sfU\cap\sfV}
\newcommand{\relFID}[2][A]{#2\operatorname{-FID}(#1)}
\newcommand{\relFPD}[2][A]{#2\operatorname{-FPD}(#1)}
\newcommand{\relsilp}[2][A]{#2\operatorname{-silp}(#1)}
\newcommand{\relspli}[2][A]{#2\operatorname{-spli}(#1)}
\newcommand{\Vid}[1][A]{\mathsf{V}\operatorname{-id}_{#1}}
\newcommand{\Upd}[1][A]{\mathsf{U}\operatorname{-pd}_{#1}}
\newcommand{\catGL}[1][\sfV]{\mathsf{LGor}_{#1}}
\newcommand{\catGR}[1][\sfU]{\mathsf{RGor}_{#1}}

\newcommand{\semiUV}{semi-$\sfU$-$\sfV$}
\newcommand{\RGorpd}[3][A]{\mathsf{RGor}_{#2}\operatorname{-pd}_{#1}#3}
\newcommand{\LGorid}[3][A]{\mathsf{LGor}_{#2}\operatorname{-id}_{#1}#3}
\newcommand{\RGorgldim}[2][A]{\mathsf{RGor}_{#2}\operatorname{-gldim}(#1)}
\newcommand{\LGorgldim}[2][A]{\mathsf{LGor}_{#2}\operatorname{-gldim}(#1)}
\newcommand{\Gfpid}[2][A]{\operatorname{Gfp-id}_{#1}#2}
\newcommand{\Gfpgldim}[1][A]{\operatorname{Gfpinj-gldim}(#1)}

\newcommand{\Cy}[2][]{\operatorname{Z}_{#1}(#2)}
\newcommand{\Co}[2][]{\operatorname{C}_{#1}(#2)}
\renewcommand{\H}[2][]{\operatorname{H}_{#1}(#2)}

\begin{document}

\title{One-sided Gorenstein rings}%

\author[L.W.\ Christensen]{Lars Winther Christensen}

\address{L.W.C. \ Texas Tech University, Lubbock, TX 79409, U.S.A.}

\email{lars.w.christensen@ttu.edu}

\urladdr{http://www.math.ttu.edu/\urltilda lchriste}

\author[S.\ Estrada]{Sergio Estrada}

\address{S.E. \ Universidad de Murcia, Murcia 30100, Spain}

\email{sestrada@um.es}

\urladdr{http://webs.um.es/sestrada}

\author[L.\ Liang]{Li Liang}

\address{L.L. \ Lanzhou Jiaotong University, Lanzhou 730070, China}

\email{lliangnju@gmail.com}

\urladdr{https://sites.google.com/site/lliangnju}

\author[P.\ Thompson]{Peder Thompson}

\address{P.T. \ M\"{a}lardalen University, V\"{a}ster{\aa}s 72123,
  Sweden}

\email{peder.thompson@mdu.se}

\urladdr{https://sites.google.com/view/pederthompson}

\author[J.\ Wang]{Junpeng Wang}

\address{J.W. \ Northwest Normal University, Lanzhou 730070, China}

\email{wangjunpeng1218@163.com}

\thanks{L.W.C.\ was partly supported by Simons Foundation
  collaboration grant 962956. S.E.\ was partly supported by grant
  PID2020-113206GB-I00 funded by MCIN/AEI/10.13039/501100011033 and
  by grant 22004/PI/22 funded by Fundaci\'on S\'eneca. L.L.\ was
  partly supported by NSF of China grant 12271230.  Part of the work
  was done during a visit by P.T.\ to the University of Murcia, the
  institution's support and hospitality is much appreciated. J.W.\ was
  partly supported by NSF of China grants 12361008 and 12061061.}

\date{11 November 2023}



\begin{abstract}
  Distinctive characteristics of Iwanaga--Gorenstein rings are
  typically understood through their intrinsic symmetry. We show that
  several of those that pertain to the Gorenstein global dimensions
  carry over to the one-sided situation, even without the noetherian
  hypothesis. Our results yield new relations among homological
  invariants related to the Gorenstein property, not only Gorenstein
  global dimensions but also the suprema of projective/injective
  dimensions of injective/projective modules and finitistic
  dimensions.
\end{abstract}

\maketitle

\thispagestyle{empty}

\section*{Introduction}

\noindent
One of the most basic homological invariants attached to a ring is the
left global dimension, that is, the supremum of the projective, or
equivalently injective, dimensions of all left modules.  For
noetherian rings the left and right global dimensions agree, but in
general they may differ, and they need not even be simultaneously
finite; see Jategaonkar \cite{AVJ69}. There are though, beyond the
noetherian realm, rings for which finiteness of the left global
dimension implies finiteness of the right global dimension. Two
results in this direction were achieved by Jensen \cite{CUJ66} and
Osofsky \cite{BLO68a}, and one of the several diverse aims of this
paper is to establish their counterparts in Gorenstein homological
algebra.

A Gorenstein analogue of rings of finite left global dimension are
those of finite left Gorenstein global dimension; Beligiannis
\cite{ABl00} calls such rings left Gorenstein, and shows that, indeed,
for a noetherian ring the left and right Gorenstein global dimensions
agree. In other words, a noetherian ring of finite left Gorenstein global
dimension is an Iwanaga--Gorenstein ring, that is, a noetherian ring
of finite self-injective dimension on both sides.  A recurring theme
of the paper---illustrated by \thmref[Theorems~]{Ggldim-new},
\thmref[]{8-UV}, \thmref[]{8-UV-dual}, \thmref[]{8}, and
\thmref[]{GwGl-1}---is that several phenomena associated with
Iwanaga--Gorenstein rings carry over to the one-sided situation even
without the noetherian hypothesis.

Let $A$ be a unital associative ring. In this paper an $A$-module
means a left $A$-module; right $A$-modules are modules over the
opposite ring $\Aop$. The Gorenstein global dimension of $A$ is
defined in terms of the Gorenstein projective dimension:
\begin{equation*}
  \Ggldim=\sup\setof{\Gpd M}{M \text{ is an $A$-module}}\:.
\end{equation*}
As in the absolute case, one can equivalently compute $\Ggldim$ based
on the Gorenstein injective dimension; see Bennis and Mahdou
\cite{DBnNMh10}. Analogous to the weak global dimension, defined in
terms of the flat dimension, the Gorenstein weak global dimension is
defined in terms of the Gorenstein flat dimension:
\begin{equation*}
  \Gwgldim = \sup\setof{\Gfd M}{M \text{ is an $A$-module}}\:.
\end{equation*}
Unlike the elementary fact that the weak global dimension of $A$ is at
most the global dimension of $A$, the corresponding inequality in
Gorenstein homological algebra,
\begin{equation*}
  \Gwgldim \leq \Ggldim \:,
\end{equation*}
was only known to hold under extra assumptions on the ring; here we
obtain it in \thmref{GwGgl_ineq}. Just as the global and weak global
dimensions, by a result of Auslander \cite{MAs55}, agree for left
noetherian rings and, trivially of course, for left perfect rings, we
observe in \corref[Corollaries~]{Bouchiba} and \corref[]{perfect} that
the equality
\begin{equation*}
  \Gwgldim = \Ggldim
\end{equation*}
holds if $A$ is left noetherian or left perfect.

While the Gorenstein weak global dimension is symmetric, see for
example \cite{CET-21a}, the Gorenstein global dimension is not. But
the fact from \cite{ABl00} that the invariants $\Ggldim$ and
$\Ggldim[\Aop]$ are simultaneously finite if $A$ is noetherian is in
\corref{GwGl-2} extended to $\aleph_0$-noetherian rings and in
\thmref{Osofsky} to rings of cardinality $\aleph_n$ for some integer
$n\ge 0$. Indeed, as in the absolute cases \cite{CUJ66,BLO68a} we show
that the invariants $\Ggldim$ and $\Ggldim[\Aop]$ are simultaneously
finite and provide bounds on their difference. Namely, if $A$ is
$\aleph_0$-noetherian, then one has:
\begin{equation*}
  |\Ggldim - \Ggldim[\Aop]| \le 1 \:,
\end{equation*}
and if $A$ has cardinality $\aleph_n$ for some integer $n \ge 0$, then one has:
\begin{equation*}
  |\Ggldim - \Ggldim[\Aop]| \le n+1 \:.
\end{equation*}

Recall that the invariant $\spli$ is
$\sup\setof{\pd{I}}{I \text{ is an injective $A$-module}}$ while
$\silp$ is defined dually. The finitistic projective dimension $\FPD$
is the supremum of projective dimensions of all $A$-modules of finite
projective dimension. The finitistic injective dimension, $\FID$, and
the finitistic flat dimension, $\FFD$, are defined similarly. It is
known from work of Beligiannis and Reiten \cite{ABlIRt07} that
$\Ggldim$ is the maximum of $\spli$ and $\silp$ and that the
inequalities $\FPD \le \silp$ and $\FID \le \spli$ hold. Thus the
equalities,
\begin{equation*}
  \max\{\spli,\FPD\}
  = \Ggldim = \max\{\silp,\FID\} \:,
\end{equation*}
which follow from \thmref{Ggldim-new}, constitute an improvement: The
question whether finiteness of $\silp$ or $\spli$ implies finiteness
of $\Ggldim$ is for an Artin algebra $A$ equivalent to the Gorenstein
Symmetry Conjecture, see \cite[Note after VII.2.7]{ABlIRt07}.

Another strain of results in the paper simplify computations of the
Gorenstein global dimensions. Under the assumption that $A$ is left
noetherian, \corref{gwgldim}, \prpref{9-1}, and \corref{Bouchiba} combine to yield
\begin{align*}
  \Ggldim 
  & = \sup\setof{\Gid M}{M \text{ is a cotorsion $A$-module}} \:.
\end{align*}
As recorded in \corref{IG} this means that if $A$ is noetherian, then
it is Iwanaga--Gorenstein if and only if every cotorsion $A$-module
has finite Gorenstein injective dimension.

Under the assumption that $A$ is left coherent, we also provide
alternative ways to compute the Gorenstein weak global dimension. In
\thmref{10} we establish a parallel in Gorenstein homological algebra
to a result of Stenstr\"om \cite{BSt70} on the weak global dimension
of left coherent rings:
\begin{align*}
  \Gwgldim  & = \sup\setof{\Gfpid M}{M~\text{is an $A$-module}} \\
            & = \sup\setof{\Gpd M}{M~\text{is a finitely presented $A$-module}} \: .           
\end{align*}
Here $\Gfpid$ denotes a Gorenstein dimension defined based on
fp-injective modules; see \dfnref{WGor}. The proof of the equalities
above proceeds in two steps: (1)~The suprema are shown to equal the
global invariant based on the Gorenstein flat-cotorsion dimension from
\cite{CELTWY-21}; (2)~this invariant is shown to agree with the
Gorenstein weak global dimension for left coherent rings. In fact,
this new global invariant lurks in the background of the proofs of all
the results discussed hitherto.
\begin{equation*}
  \ast \ \ \ast \ \ \ast
\end{equation*}
On the organization of the paper: The results we have advertised above
all derive from four separate specializations of one underlying
theory. It deals with the abstract notions of Gorenstein modules
associated to a cotorsion pair as first defined in \cite{CET-20}, and
the associated relative homological dimensions developed in
\cite{CDEHLT-23,CELTWY-21}.

In \secref[Sections~]{uv} and \secref[]{relgor} we work in the setting
of a hereditary cotorsion pair $\sfUV$ in $\Mod$, which is generated
by a set and exhibits periodicity, see \dfnref{periodic}. The
motivation comes from the cotorsion pair $(\Flat,\Cot)$ and an
associated homological invariant, the Gorenstein flat-cotorsion
dimension we mentioned above.

The main results of the paper are derived in
\secref[Sections~]{gfcgldim}--\secref[]{gwgldim} by specializing the
results of the first two sections to the cotorsion pairs
\begin{equation*}
  (\Prj,\Mod) \,, \quad\, (\Mod,\Inj)\,,\quad (\Flat,\Cot)\,,
\end{equation*}
and, for left coherent rings, to the pair
$({}^\perp\FpInj,\FpInj)$. The last specialization produces the result
inspired by Stenstr\"om that was mentioned above, but the majority of
the results come from the specialization to $(\Flat,\Cot)$, which also
contributes to the general theory of the Gorenstein flat-cotorsion
dimension.


\section{Relative homological dimensions under periodicity}
\label{sec:uv}

\noindent
Throughout the paper, $A$ denotes an associative unital ring.  By an
$A$-module we mean a left $A$-module; right $A$-modules are modules
over the opposite ring $\Aop$. The category of $A$-modules is denoted
$\Mod$. Chain complexes of $A$-modules, or $A$-complexes, are indexed
homologically. That is, the degree $n$ component of an $A$-complex $M$
is denoted $M_n$, and we denote the degree $n$ cycle, cokernel, and
homology modules by $\Cy[n]{M}$, $\Co[n]{M}$, and $\H[n]{M}$,
respectively.

For a class $\sfC$ of $A$-modules we denote by $\sfC^\perp$ the class
of $A$-modules $M$ with $\Ext^1_A(C,M) = 0$ for all $C \in \sfC$;
similarly, the class ${}^\perp\sfC$ consists of the $A$-modules $M$
with $\Ext^1_A(M,C) = 0$ for all $C \in \sfC$.  Recall that a pair
$\sfUV$ of classes of $A$-modules constitutes a \emph{cotorsion pair}
if $\sfU = {}^\perp\sfV$ and $\sfU^\perp = \sfV$ hold. A cotorsion
pair $\sfUV$ is \emph{complete} if every $A$-module $M$ has a special
$\sfU$-precover and a special $\sfV$-preenvelope, i.e.\ there are
exact sequences
\begin{equation*}
  0 \lra V \lra U \lra M \lra 0 \qquad\text{and}\qquad
  0 \lra M \lra V' \lra U' \lra 0
\end{equation*}
with $U,U'$ in $\sfU$ and $V,V'$ in $\sfV$. A cotorsion pair $\sfUV$
is \emph{generated by a set} if there is a set $\sfS$ of $A$-modules
with $\sfS^\perp=\sfV$.  By a result of Eklof and Trlifaj
\thmcite[10]{PCEJTr01}, a cotorsion pair that is generated by a set is
complete.  A cotorsion pair $\sfUV$ is \emph{hereditary} if
$\Ext_A^i(U,V)=0$ holds for all $U$ in $\sfU$, all $V$ in $\sfV$, and
all integers $i \ge 1$. Equivalently, $\sfUV$ is hereditary if and
only if the class $\sfU$ is resolving (i.e.\ closed under kernels of
surjective homomorphisms) and if and only if $\sfV$ is coresolving (i.e.\ closed
under cokernels of injective homomorphisms).

Throughout this section, $\sfUV$ is a hereditary cotorsion pair in
$\Mod$, and we assume that it is generated by a set.  As in
\cite{CDEHLT-23} we say that an $A$-complex $V$ is
\emph{$\sfV$-acyclic} if it is acyclic and $\Cy[n]{V}$ belongs to
$\sfV$ for all $n\in \ZZ$, and an $A$-complex $U$ of modules in $\sfU$
is \emph{semi-$\sfU$} if $\Hom[A](U,V)$ is acyclic for every
$\sfV$-acyclic complex $V$. Similarly one defines $\sfU$-acyclic
complexes and semi-$\sfV$ complexes.

\begin{dfn}
  \label{dfn:periodic}
  We say that $\sfUV$ is \emph{right periodic} if every acyclic
  complex of modules from $\sfV$ is $\sfV$-acyclic and \emph{weakly
    right periodic} if every acyclic complex of injective $A$-modules
  is $\sfV$-acyclic. Dually, we say that $\sfUV$ is \emph{left
    periodic} if every acyclic complex of modules from $\sfU$ is
  $\sfU$-acyclic and \emph{weakly left periodic} if every acyclic
  complex of projective $A$-modules is $\sfU$-acyclic.
\end{dfn}

Notice that if $\sfUV$ is right/left periodic, then it is weakly
right/left periodic.  The terminology in the definition above comes
from the fact that $\sfUV$ is \rightbothperiodic\ if and only if in
any exact sequence of $A$-modules of the form
\begin{equation*}
  0 \lra M \lra V \lra M \lra 0
\end{equation*}
with $V$ in $\sfV$ (or injective) the module $M$ belongs to
$\sfV$. The \leftbothperiodic\ property can be described similarly.
See for example Bazzoni, Cort\'es Izurdiaga, and Estrada
\prpcite[2.4]{BCE-20}. Examples of left/right periodic cotorsion pairs
come up at the end of this section and in the opening paragraphs of
\secref[Sections~]{gfcgldim} and \secref[]{stenstrom}.

Every $A$-complex $M$ admits a \emph{semi-$\sfU$ replacement}, that
is, a semi-$\sfU$ complex $U$ that is isomorphic to $M$ in the derived
category over $A$. Every $A$-complex also admits a similarly defined
\emph{semi-$\sfV$ replacement.} As a matter of fact, one can even get
resolutions; that is, quasi-isomorphisms $U \qra M$ and $M \qra V$ in
the category of $A$-complexes where $U$ and $V$ are semi-$\sfU$ and
semi-$\sfV$ complexes. See for example Gillespie \cite{JGl04} and Yang
and Liu \thmcite[3.5]{GYnZLi11}. For an $A$-module $M$ one can
construct a semi-$\sfU$ replacement by taking successive special
$\sfU$-precovers and a semi-$\sfV$ replacement by taking successive
special $\sfV$-preenvelopes in which case one gets
resolutions $U \qra M$ and $M \qra V$ in the category of
$A$-complexes. An $A$-complex that is both
semi-$\sfU$ and semi-$\sfV$ is called \emph{semi-}$\sfU$-$\sfV$. A
\emph{semi-$\sfU$-$\sfV$ replacement} of an $A$-complex $M$ is a
semi-$\sfU$-$\sfV$ complex that is isomorphic to $M$ in the derived
category over $A$. See \seccite[4]{CDEHLT-23}.

We first handle some aspects of the relative homological dimension
with respect to the class $\sfU$, primarily in the case where $\sfUV$
is \wrightperiodic, and then from \lemref{Vn} we shift focus to $\sfV$
and the weakly left periodic property.

\begin{dfn}
  Let $M \ne 0$ be an $A$-module. The \emph{$\sfU$-projective
    dimension} of $M$, denoted $\Upd M$, is the least integer $n$ such
  that there exists a semi-$\sfU$ resolution $U \qra M$ with $U_i=0$
  for $i>n$.  If no such resolution exists, then one sets
  $\Upd M=\infty$; further, one sets $\Upd 0 = -\infty$.  We say that
  a module $M$ has \emph{finite} $\sfU$-projective dimension if
  $\Upd{M} < \infty$ holds.
\end{dfn}

We denote by $\Prj$ the class of projective $A$-modules. For the
cotorsion pair $(\sfU,\sfV) = (\Prj,\Mod)$ we of course use the
standard notation $\mathrm{pd}_A$ for the $\sfU$-projective dimension.
The equivalence of some of the assertions in the next result is
standard; we include an entire proof for completeness.

\begin{prp}
  \label{prp:Upd}
  Let $M$ be an $A$-module and $n \ge 0$ an integer. The following
  conditions are equivalent.
  \begin{eqc}
  \item $\Upd M \leq n$.
  \item There exists a semi-$\sfU$-$\sfV$ replacement $W$ of $M$ with
    $W_i=0$ for $i>n$.
  \item For every semi-$\sfU$-$\sfV$ replacement $W$ of $M$, the
    cokernel $\Co[n]{W}$ is in $\capUV$.
  \item For every semi-$\sfU$ replacement $U$ of $M$, the cokernel
    $\Co[n]{U}$ is in $\sfU$.
  \item For every projective resolution $P \qra M$, the cokernel
    $\Co[n]{P}$ is in $\sfU$.
  \item For every $V$ in $\sfV$, one has $\Ext_A^{n+1}(M,V)=0$.
  \end{eqc}
  Moreover, there is an equality,
  \begin{equation*}
    \Upd{M} = \sup\setof{i \in \ZZ}{\Ext_A^i(M,V) \ne 0 \text{ for some module }
      V \in \sfV} \:,
  \end{equation*}
  and if two out of three modules in a short exact sequence have
  finite $\sfU$-projective dimension, then so has the third.
\end{prp}

\begin{prf*}
  \proofofimp{i}{ii} Let $U \qra M$ be a semi-$\sfU$ resolution with
  $U_i=0$ for $i>n$. By \thmcite[A.8]{CDEHLT-23} there is an exact
  sequence $0\to U \to W\to U'\to 0$ of $A$-complexes such that $U'$
  is $\sfU$-acyclic and $W$ is semi-$\sfU$-$\sfV$ with $W_i=0$ for
  $i\gg0$. Acyclicity of $U'$ yields an exact sequence
  $0\to \Co[n]{U}\to \Co[n]{W}\to \Co[n]{U'}\to 0$.  As the outer
  terms belong to $\sfU$, so does $\Co[n]{W}$. As $\sfUV$ is
  hereditary and $W_i=0$ holds for $i\gg0$ it follows that that
  $\Co[i]{W}$ is in fact in \mbox{$\capUV$} for $i \ge n$. It now
  follows that the exact sequence
  $0 \to \Co[n+1]{W}\to W_n\to \Co[n]{W}\to 0$ splits and, therefore,
  we may replace $W$ by the homotopic semi-$\sfU$-$\sfV$ complex
  $0 \to \Co[n]{W}\to W_{n-1}\to \cdots$, which is the desired
  \semiUV\ replacement of $M$.
  
  \proofofimp{ii}{iii} Given any other semi-$\sfU$-$\sfV$ replacement
  $W'$ of $M$, \prpcite[A.7]{CDEHLT-23} says that $\Co[n]{W'}$ is a
  direct summand of $\Co[n]{W} \oplus X = W_n \oplus X$ where $X$ is a
  module in $\capUV$. Thus $\Co[n]{W'}$ belongs to $\capUV$.
  
  \proofofimp{iii}{iv} Let $U$ be any semi-$\sfU$ replacement of
  $M$. Again by \thmcite[A.8]{CDEHLT-23} there is an exact sequence of
  $A$-complexes $0 \to U\to W \to U'\to 0$ where $W$ is
  semi-$\sfU$-$\sfV$ and $U'$ is $\sfU$-acyclic, and again acyclicity
  of $U'$ yields an exact sequence
  $0\to \Co[n]{U}\to \Co[n]{W}\to \Co[n]{U'}\to 0$. The middle and
  right-hand modules belong to $\sfU$, and thus the resolving property
  of $\sfU$ implies that $\Co[n]{U}$ belongs to $\sfU$.
  
  \proofofimp{iv}{v} This is clear since projective resolutions yield
  semi-$\sfU$ replacements.
  
  \proofofimp{v}{vi} Let $V$ be in $\sfV$, and take a projective
  resolution $P \qra M$. Dimension shifting yields
  $\Ext_A^{n+1}(M,V)\cong \Ext_A^1(\Co[n]{P},V)=0$, since $\Co[n]{P}$
  belongs to $\sfU$.
  
  \proofofimp{vi}{i} Construct a semi-$\sfU$ resolution $U \qra M$
  from special $\sfU$-precovers. For $V$ in $\sfV$, one has
  $0=\Ext_A^{n+1}(M,V)\cong \Ext_A^1(\Co[n]{U},V)$, thus $\Co[n]{U}$
  is in $\sfU$.

  The equality and the last assertion are immediate from the
  equivalence of conditions \eqclbl{i} and \eqclbl{vi}.
\end{prf*}

\begin{dfn}
  Let $n \ge 0$ be an integer. Denote by $\sfU_n$ the class of
  $A$-modules of $\sfU$-projective dimension at most $n$.
\end{dfn}

\begin{fct}
  \label{fct:Un}
  Let $n \ge 0$ be an integer. The pair $(\sfU_n,\sfU_n^\perp)$ is by
  work of Cort\'es Izurdiaga, Estrada, and Guil Asensio \cite[Thm.\
  2.2, Cor.\ 2.3 and 2.7]{CEG-12} a hereditary cotorsion pair
  generated by a set; in particular, it is complete.
\end{fct}

The next proposition shows that \rightbothperiodic\ cotorsion pairs
come in families.
\begin{prp}
  \label{prp:periodic-n}
  Let $n\ge 0$ be an integer. If $\sfUV$ is \rightbothperiodic, then
  $(\sfU_n,\sfU_n^\perp)$ is \rightbothperiodic. In particular, if
  $\sfUV$ is \wrightperiodic, then every Gorenstein injective
  $A$-module belongs to $\sfU_n^\perp$.
\end{prp}

\begin{prf*}
  First consider the case where $\sfUV$ is \wrightperiodic. Let $I$ be
  an acyclic complex of injective $A$-modules and $L$ a module in
  $\sfU_n$. For any module $V$ in $\sfV$, note that \prpref{Upd}
  yields $\Ext_A^{n+1}(L,V)=0$. Now, for every $i\in \ZZ$ there is an
  exact sequence
  \begin{equation*}
    0 \lra \Cy[i+n]{I} \lra I_{i+n}\lra \cdots \lra I_{i+1}\lra
    \Cy[i]{I}\lra 0 \:,
  \end{equation*}
  so dimension shifting yields
  $\Ext_A^1(L,\Cy[i]{I})\cong \Ext_A^{n+1}(L,\Cy[i+n]{I})=0$, as
  $\Cy[i+n]{I}$ belongs to $\sfV$ by hypothesis. Thus $\Cy[i]{I}$
  belongs to $\sfU_n^{\perp}$. The in particular statement follows
  simply because every Gorenstein injective $A$-module is a cycle
  submodule in an acyclic complex of injective $A$-modules.

  The case where $\sfUV$ is \rightperiodic\ is handled in the same
  way, only one takes $I$ to be an acyclic complex of modules from
  $\sfU_n^\perp$.
\end{prf*}

\begin{dfn}
  Imitating the finitistic projective dimension, set
  \begin{align*}
    \relFPD{\sfU} &= \sup\setof{\Upd M}
                    {M \text{ is an $A$-module with } \Upd M<\infty} \:.
  \end{align*}
  Further, imitating the invariant $\spli$, set
  \begin{align*}
    \relspli{\sfU}&=\sup\setof{\Upd I}
                    {I \text{ is an injective $A$-module}} \:.
  \end{align*}
\end{dfn}

The next lemma is a more general version of Emmanouil's
\lemcite[5.2]{IEm12}.

\begin{lem}
  \label{lem:Emm}
  Assume that $\relspli{\sfU} = n < \infty$ holds and let $M$ be an
  $A$-module.  For every projective resolution $P \qra M$ there is an
  acyclic complex $U$ of modules in $\sfU$ with $U_i=P_i$ for
  $i\geq n$.
\end{lem}

\begin{prf*}
  For the class $\sfU$ of flat $A$-modules, the statement is proved in
  \lemcite[5.2]{IEm12}. A straightforward generalization of the
  argument provided there yields the claim asserted here.
\end{prf*}

In the next statement, and in the rest of the paper, $\Inj$ and $\GI$
denote the classes of injective and Gorenstein injective
$A$-modules. The purpose of the statement is to tie in the abstract
notion of modules of finite $\sfU$-projective dimension with the
better known notions of injective and Gorenstein injective modules.

\begin{prp}
  \label{prp:4-UV}
  Let $n \ge 0$ be an integer and assume that $\sfUV$ is
  \wrightperiodic. The following conditions are equivalent.
  \begin{eqc}
  \item $\max\{\relspli{\sfU},\relFPD{\sfU}\} \leq n$.
  \item $\relspli{\sfU} = \relFPD{\sfU} \leq n$.
  \item ${\sfU}_n \cap {\sfU}_n^\perp = \Inj$.
  \item ${\sfU}_n^\perp = \GI$.
  \end{eqc}
\end{prp}

\begin{prf*}
  \proofofimp{i}{ii} Assuming that both $\relspli{\sfU}$ and
  $\relFPD{\sfU}$ are finite, the inequality
  $\relspli{\sfU} \le \relFPD{\sfU}$ is evident. For the opposite
  inequality, set $s = \relspli{\sfU}$ and consider a module $M$ of
  finite $\sfU$-projective dimension. It follows from \lemref{Emm}
  that the $s^\mathrm{th}$ cokernel module in a projective resolution
  of $M$ appears as a cokernel module in an acyclic complex of modules
  from $\sfU$. Finiteness of $\relFPD{\sfU}$ implies that this module
  belongs to $\sfU$: Indeed, it is for every $m \ge 0$ an
  $m^\mathrm{th}$ syzygy of a module of finite $\sfU$-projective
  dimension.

  \proofofimp{ii}{iii} The inclusion ``$\supseteq$'' holds as
  injective $A$-modules by assumption have $\sfU$-projective dimension
  at most $n$.  For the opposite inclusion, let $M$ be a module in the
  intersection ${\sfU}_n \cap {\sfU}_n^\perp$.  Consider an exact
  sequence $0 \to M \to I \to N \to 0$ of $A$-modules with $I$
  injective.  Since $I$ and $M$ have finite $\sfU$-projective
  dimension, so has $N$; see \prpref{Upd}. It follows from the
  assumption $\relFPD{\sfU} \leq n$ that $N$ is in $\sfU_n$, and so
  one has $\Ext_A^{1}(N,M)=0$.  Thus the sequence
  $0 \to M \to I \to N \to 0$ is split exact, whence $M$ is injective.

  \proofofimp{iii}{iv} Per \prpref{periodic-n} one has
  ${\sfU}_n^{\perp} \supseteq \GI$. For the other containment, let $M$
  belong to ${\sfU}_n^{\perp}$. By \fctref{Un},
  $(\sfU_n , \sfU_n^\perp)$ is a complete cotorsion pair, so there is
  an exact sequence $0 \to M' \to U \to M \to 0$ with $U \in \sfU_n$
  and $M' \in \sfU_n^\perp$. Thus $U$ also belongs to $\sfU_n^\perp$,
  and hence it is injective by hypothesis. Repeating this process and
  taking an injective resolution of $M$ yields an acyclic complex of
  injective modules with $M$ as a cycle module. Any such complex $I$
  is totally acyclic: Indeed, the cycle modules in such a complex
  belong to the class $\sfV$ as $\sfUV$ is \wrightperiodic, and by
  assumption injective modules belong to $\sfU_n$, so for an injective
  module $E$ one has
  $\Ext_A^1(E,\Cy[i]{I}) \cong \Ext_A^{n+1}(E,\Cy[i+n]{I}) = 0$ by
  \prpref{Upd}.

  \proofofimp{iv}{i} As $(\sfU_n,\GI)$ per the assumption and
  \fctref{Un} is a cotorsion pair, one has
  \begin{equation}
    \tag{$\ast$}
    {}^\perp\GI = {}^\perp(\sfU_n^\perp) = \sfU_n \:.
  \end{equation}
  As injective modules belong to ${}^\perp\GI$ one has
  $\relspli{\sfU} \leq n$. To see that $\relFPD{\sfU} \leq n$ holds,
  let $M$ be an $A$-module of finite $\sfU$-projective dimension.  It
  follows from \prpref{periodic-n} that $\Ext_A^{1}(M, G)=0$ holds for
  every $G \in \GI$, whence $M$ is in ${\sfU}_n$ by $(\ast)$.
\end{prf*}

\begin{prp}
  \label{prp:7-UV}
  Assume that $\sfUV$ is \wrightperiodic\ and let $V$ be a module in
  $\sfV$. If\, $\Upd{V}$ is finite, then the following inequality
  holds:
  \begin{equation*}
    \id V \leq \max\{\relspli{\sfU}, \relFPD{\sfU} \} \:.
  \end{equation*}
\end{prp}

\begin{prf*}
  We may assume that $\max\{\relspli{\sfU}, \relFPD{\sfU} \}=n$ holds
  for some $n \ge 0$.  Consider an exact sequence of $A$-modules,
  \begin{equation*}
    0 \lra V \lra I_{0} \lra I_{-1} \lra \cdots \lra I_{-n+1} \lra G\ \lra 0
  \end{equation*}
  where each module $I_{i}$ is injective. For every module $M$ of
  finite $\sfU$-projective dimension $\Ext_A^{n+1}(M,V)=0$ holds by
  \prpref{Upd}, and dimension shifting yields
  $\Ext_A^{1}(M,G)\cong \Ext_A^{n+1}(M,V)=0$, so that $G$ belongs to
  $\sfU_n^\perp$.  Since the modules $V,I_0,\ldots,I_{-n+1}$ have
  finite $\sfU$-projective dimension, so does $G$, see
  \prpref{Upd}. Thus $G$ belongs to \mbox{$\sfU_n \cap \sfU_n^\perp$}
  and is, therefore, injective by \prpref{4-UV}.
\end{prf*}

Dual to the $\sfU$-projective dimension one defines the
\emph{$\sfV$-injective dimension} of an $A$-module $M$; the notation
is $\Vid M$ and for an integer $n \ge 0$ the class of $A$-modules of
$\sfV$-injective dimension at most $n$ is denoted by $\sfV_n$.  For
the cotorsion pair $(\sfU,\sfV) = (\Mod,\Inj)$ we of course use the
standard notation $\mathrm{id}_A$, as already done in \prpref{7-UV},
for the $\sfV$-injective dimension.
  
The $\sfV$-injective dimension has properties dual to those in
\prpref{Upd}:

\begin{prp}
  \label{prp:Vid}
  Let $M$ be an $A$-module and $n \ge 0$ an integer. The following
  conditions are equivalent.
  \begin{eqc}
  \item $\Vid M \leq n$.
  \item There exists a semi-$\sfU$-$\sfV$ replacement $W$ of $M$ with
    $W_{-i}=0$ for $i>n$.
  \item For every semi-$\sfU$-$\sfV$ replacement $W$ of $M$, the
    kernel $\Cy[-n]{W}$ is in $\capUV$.
  \item For every semi-$\sfV$ replacement $V$ of $M$, the kernel
    $\Cy[-n]{V}$ is in $\sfV$.
  \item For every injective resolution $M \qra I$, the kernel
    $\Cy[-n]{I}$ is in $\sfV$.
  \item For every $U$ in $\sfU$, one has $\Ext_A^{n+1}(U,M)=0$.
  \end{eqc}
  Moreover, there is an equality,
  \begin{equation*}
    \Vid{M} = \sup\setof{i \in \ZZ}
    {\Ext_A^i(U,M) \ne 0 \text{ for some module }
      U \in \sfU} \:,
  \end{equation*}
  and if two out of three modules in a short exact sequence have
  finite $\sfV$-injective dimension, then so has the third.
\end{prp}

\begin{prf*}
  Dual to the proof of \prpref{Upd}.
\end{prf*}

\begin{dfn}
  Imitating the finitistic injective dimension, set
  \begin{align*}
    \relFID{\sfV} &= \sup\setof{ \Vid M }
                    { M \text{ is an $A$-module with } \Vid M<\infty } \:.
  \end{align*}
  Further, imitating the invariant $\silp$, set
  \begin{align*}
    \relsilp{\sfV}
    &=\sup\setof{\Vid P}{P \text{ is a projective $A$-module}} \:.
  \end{align*}
\end{dfn}

To investigate the invariants introduced right above, we start by
establishing a counterpart to \fctref{Un}. The proof stands out from
the proofs of \prpref[]{Vid} and
\prpref[]{periodic-n-bis}--\prpref[]{7-UV-bis} by not being simply
dual to the proof of \fctref{Un}. Our argument below is modeled on the
proof of \thmcite[8.7]{RGbJTr12-1}.

\begin{lem}
  \label{lem:Vn}
  Let $n \ge 0$ be an integer. The pair $({}^\perp\sfV_n,\sfV_n)$ is a
  complete hereditary cotorsion pair generated by a set.
\end{lem}

\begin{prf*}
  Recall that $\sfUV$ is generated by a set $\sfS$. For each
  $S\in \sfS$, fix a projective resolution $P_S\qra S$. We argue that
  the set
  \begin{equation*}
    \sfS' = \setof{\Co[n]{P_S}}{S\in \sfS}
  \end{equation*}
  generates $({}^\perp\sfV_n,\sfV_n)$, that is,
  $\sfS'^{\perp}=\sfV_n$. Let $M$ be an $A$-module and fix an
  injective resolution $M\qra I$. Now $M\in \sfS'^\perp$ means
  $\Ext_A^1(\Co[n]{P_S},M)=0$ for every $S\in \sfS$. However, one has
  \begin{equation*}
    \Ext_A^1(\Co[n]{P_S},M) \cong \Ext_A^{n+1}(S,M)
    \cong \Ext_A^{1}(S,\Cy[-n]{I}) \:,
  \end{equation*}
  which vanishes if and only if $\Cy[-n]{I}$ belongs to $\sfV$, that
  is, if and only if $M$ belongs to $\sfV_n$, see \prpref{Vid}. This
  shows that $({}^\perp\sfV_n,\sfV_n)$ is a cotorsion pair generated
  by $\sfS'$. To see that it is hereditary, let
  $0 \to M' \to M \to M''\to 0$ be an exact sequence with
  $M',M\in \sfV_n$. Given injective resolutions $M'\qra I'$ and
  $M''\qra I''$, the Horseshoe Lemma yields an injective resolution
  $M \qra I$, such that the sequence $0 \to I' \to I \to I'' \to 0$ is
  exact. Thus one obtains an exact sequence
  $0\to \Cy[-n]{I'} \to \Cy[-n]{I} \to \Cy[-n]{I''}\to 0$. As
  $\Cy[-n]{I'}$ and $\Cy[-n]{I}$ belong to $\sfV$, so does
  $\Cy[-n]{I''}$ by the coresolving property of $\sfV$. Thus $M''$
  belongs to $\sfV_n$; that is, $\sfV_n$ is coresolving and it follows
  that the pair $({}^\perp\sfV_n,\sfV_n)$ is hereditary.  Finally, a
  hereditary cotorsion pair generated by a set is complete.
\end{prf*}

\begin{prp}
  \label{prp:periodic-n-bis}
  Let $n\ge 0$ be an integer. If $\sfUV$ is \leftbothperiodic, then
  $({}^\perp\sfV_n,\sfV_n)$ is \leftbothperiodic. In particular, if
  $\sfUV$ is \wleftperiodic, then every Gorenstein projective
  $A$-module belongs to ${}^\perp\sfV_n$.
\end{prp}

\begin{prf*}
  Dual to the proof of \prpref{periodic-n}.
\end{prf*}

\begin{lem}
  \label{lem:Emm-bis}
  Assume that $\relsilp{\sfV}=n<\infty$ holds and let $M$ be an
  $A$-module. For every injective resolution $M\qra I$ there is an
  acyclic complex $V$ of modules in $\sfV$ with $V_{-i}=I_{-i}$ for
  $i\geq n$.
\end{lem}

\begin{prf*}
  The argument in the proof of \lemcite[5.2]{IEm12} generalizes and
  dualizes to prove the claim, see also the proof of \lemref{Emm}.
\end{prf*}

In the next statement, and in the rest of the paper, $\GP$ denotes the
class of Gorenstein projective $A$-modules.

\begin{prp}
  \label{prp:4-UV-bis}
  Let $n \ge 0$ be an integer and assume that $\sfUV$ is
  \wleftperiodic. The following conditions are equivalent.
  \begin{eqc}
  \item $\max\{\relsilp{\sfV},\relFID{\sfV}\} \leq n$.
  \item $\relsilp{\sfV} = \relFID{\sfV} \leq n$.
  \item ${}^\perp{\sfV}_n \cap {\sfV}_n = \Prj$.
  \item ${}^\perp{\sfV}_n = \GP$.
  \end{eqc}
\end{prp}

\begin{prf*}
  Dual to the proof of \prpref{4-UV}.
\end{prf*}

\begin{prp}
  \label{prp:7-UV-bis}
  Assume that $\sfUV$ is \wleftperiodic\ and let $U$ be a module in
  $\sfU$. If\, $\Vid{U}$ is finite, then the following inequality
  holds:
  \begin{equation*}
    \pd U \leq \max\{\relsilp{\sfV}, \relFID{\sfV} \} \:.
  \end{equation*}
\end{prp}

\begin{prf*}
  Dual to the proof of \prpref{7-UV}.
\end{prf*}

We close this section with an application of the theory developed up
to now, more specifically, we apply \prpref[Propositions~]{7-UV} and
\prpref[]{7-UV-bis} to the absolute cotorsion pairs $(\Prj,\Mod)$ and
$(\Mod,\Inj)$. It is trivial that the former is \rightperiodic\ and
the latter is \leftperiodic. For the former, the relevant special case
of \fctref{Un} is an earlier result of Aldrich, Enochs, Jenda, and
Oyonarte \thmcite[4.2]{AEJO-01}.

It is known, e.g.\ from Beligiannis and Reiten \prpcite[VII.1.3 and
Thm.\ VII.2.2]{ABlIRt07} or Emmanouil \thmcite[4.1]{IEm12}, that
$\Ggldim$ is finite if and only if the invariants $\spli$ and $\silp$
both are finite. As the inequalities $\FPD \le \silp$ and
$\FID \le \spli$ always hold, also by \prpcite[VII.1.3]{ABlIRt07}, the
next theorem is an improvement on this characterization of rings of
finite Gorenstein global dimension.  To parse the full statement
recall the invariant
\begin{equation*}
  \sfli = \sup\setof{\fd{I}}{I \text{ is an injective $A$-module}} \:.
\end{equation*}
Here, of course, $\mathrm{fd}_A$ is the standard notation for the flat
dimension.

\begin{thm}
  \label{thm:Ggldim-new}
  Let $n \ge 0$ be an integer. The following conditions are
  equivalent.
  \begin{eqc}
  \item $\Ggldim\leq n$.
  \item $\max\{\spli,\FPD\}\leq n$.
  \item $\max\{\silp,\FID\}\leq n$.
  \item $\max\{\sfli,\FPD\}\leq n$.
  \end{eqc}
\end{thm}

\begin{prf*}
  It is known from \prpcite[VII.1.3 and Thm.\ VII.2.2]{ABlIRt07} and
  \thmcite[4.1]{IEm12} that \eqclbl{i} implies both \eqclbl{ii} and
  \eqclbl{iii}. It follows from \prpref{7-UV} applied to the cotorsion
  pair $(\Prj,\Mod)$ that under the assumptions in \eqclbl{ii} also
  $\silp \le n$ holds, which by the references above means that
  $\Ggldim$ is at most $n$. Similarly, it follows from
  \prpref{7-UV-bis} applied to the cotorsion pair $(\Mod,\Inj)$ that
  under the assumptions in \eqclbl{iii} also $\spli \le n$ holds, and
  then $\Ggldim \le n$ holds by the references above. This shows the
  equivalence of \eqclbl{i}--\eqclbl{iii}.
  
  Condition \eqclbl{ii} implies \eqclbl{iv} simply because
  $\sfli\le \spli$ holds. Conversely, notice that $\sfli < \infty$ and $\FPD<\infty$
  imply $\spli < \infty$ by a result of
  Emmanouil and Talelli \prpcite[2.2]{IEmOTl11}. Thus one has
  $\spli\le \FPD$.
\end{prf*}

\begin{rmk}
  One could add to the theorem above another equivalent condition,
  \begin{eqc}\setcounter{eqc}{4}
  \item $\max\{\silf,\FID\}\leq n$,
  \end{eqc}
  but it coincides with \eqclbl{iii} as $\silf=\silp$ holds by
  \prpcite[2.1]{IEmOTl11}.
\end{rmk}


\section{Relative Gorenstein global dimensions}
\label{sec:relgor}

\noindent
As in the previous section, $\sfUV$ is a hereditary cotorsion pair in
$\Mod$ generated by a set, in particular, it is complete.

Recall from \dfncite[2.1]{CET-20} that an $A$-module $M$ is called
\emph{right $\sfU$-Gorenstein} if one has $M=\Co[0]{T}$ for some right
$\sfU$-totally acyclic complex $T$. That is to say, $T$ is an acyclic
complex of modules in $\sfU$ with cycle modules in $\sfV$, such that
$\Hom(T,W)$ is acyclic for every module $W$ in $\capUV$. Denote by
$\catGR$ the category of right $\sfU$-Gorenstein $A$-modules.  We
remind the reader that the $\catGR$-projective dimension of an
$A$-module $M\ne 0$, denoted $\RGorpd{\sfU}{M}$, was defined in
\dfncite[4.5]{CDEHLT-23} as the least integer $n \ge 0$ such that
there is a semi-$\sfU$-$\sfV$ replacement $W$ of $M$ with $\Co[n]{W}$
in $\catGR$. If no such integer exists, then the $\catGR$-projective
dimension of $M$ is infinite, and one sets
$\RGorpd{\sfU}{0} = -\infty$. We say that a module $M$ has
\emph{finite} $\catGR$-projective dimension if
$\RGorpd{\sfU}{M} < \infty$ holds.

For $(\sfU,\sfV) = (\Prj,\Mod)$ the $\catGR$-projective dimension is
the standard Gorenstein projective dimension and, of course, we use
the standard notation $\mathrm{Gpd}_A$.

\begin{lem}
  \label{lem:dimineq}
  Let $M$ be an $A$-module. There is an inequality,
  \begin{equation*}
    \RGorpd{\sfU}{M} \leq \Upd M \:,
  \end{equation*}
  and equality holds if $\Upd M$ is finite.
\end{lem}

\begin{prf*}
  We may assume that $M$ is nonzero.  Assume that $\Upd M\leq n$ holds
  for an integer $n$. There is by \prpref{Upd} a \semiUV\ replacement
  $W \qis M$ with $W_i = 0$ for $i > n$.  As $\Co[n]{W} = W_n$ belongs
  to $\capUV$ and hence to $\catGR$, see \exacite[2.2]{CET-20}, one
  has $\RGorpd{\sfU}{M}\leq n$.

  To see that equality holds, notice that if $\Co[m]{W}$ belongs to
  $\catGR$ for some $m < n$, then one has an exact sequence,
  \begin{equation*}
    0 \lra W_{n} \lra \cdots \lra W_m \lra \Co[m]{W} \lra 0 \:,
  \end{equation*}
  which breaks into split short exact sequences, as
  $\Ext^i_A(\Co[m]{W},W_n) = 0$ holds for all $i \ge 1$ by
  \rmkcite[4.9]{CDEHLT-23}. Thus also $\Co[m]{W}$ belongs to $\sfU$.
\end{prf*}

\begin{thm}
  \label{thm:SymU}
  Let $M$ be an $A$-module. If\, $\id M$ is finite, then one has
  \begin{equation*}
    \RGorpd{\sfU}{M} = \Upd M \:.
  \end{equation*}
\end{thm}
The second paragraph below follows the proof of \thmcite[2.1]{CET-21a}
quite closely. We have included the full argument as we will later
direct the reader to dualize it.

\begin{prf*}
  We may assume that $M$ is nonzero.  Set $n=\id M$ and let $M\qra I$
  be an injective resolution with $I_{-i}=0$ for $i>n$. As one has
  $\Vid{M} \le \id{M} = n$ hold, there is by \prpref{Vid} a
  semi-$\sfU$-$\sfV$ replacement $W \qis I$ with $W_{-i}=0$ for
  $i>n$. By Avramov and Foxby \cite[1.4.I]{LLAHBF91} there is a
  quasi-isomorphism $W \qra I$ in the category of $A$-complexes.

  Fix $m\ge \RGorpd{\sfU}{M}$; the module $\Co[m]{W}$ belongs to
  $\catGR$ by \lemcite[4.7]{CDEHLT-23}. We argue next that
  $\Ext^1_A(G,\Co[m]{W}) = 0$ holds for every $A$-module $G$ in
  $\catGR$. Fix such a module, by definition there is an exact
  sequence of $A$-modules,
  \begin{equation*}
    0 \lra G \lra U_0 \lra \cdots \lra U_{-m-n+1} \lra G' \lra 0 \:,
  \end{equation*}
  with each module $U_i$ in $\capUV$ and $G'$ in $\catGR$. As
  $\Co[m]{W}$ belongs to $\sfV$, dimension shifting yields:
  \begin{equation}
    \tag{$\ast$}
    \Ext_A^1(G,\Co[m]{W}) \cong \Ext_A^{m+n+1}(G’,\Co[m]{W}) \:.
  \end{equation}
  Let $C$ be the mapping cone of the quasi-isomorphism $W \qra I$; one
  has $C_i=0$ for $i < -n$, and the complex consists of direct sums of
  modules that are in $\capUV$ or injective. Moreover, one has
  $C_{-n} = I_{-n}$ as $W_i = 0$ for $i<-n$. Further, as $I_v=0$ holds
  for $v > 0$ and one has $m \ge 0$, there is an isomorphism
  $\Co[m+1]{C} \cong \Co[m]{W}$.  Thus dimension shifting along the
  exact sequence
  \begin{equation*}
    0 \lra \Co[m+1]{C} \lra C_{m} \lra \cdots \lra C_{-n+1} \lra I_{-n} \lra 0
  \end{equation*}
  yields
  \begin{equation}
    \tag{$\dagger$}
    \Ext_A^{m+n+1}(G’,\Co[m]{W}) \cong \Ext_A^{1}(G’,I_{-n}) = 0 \:.
  \end{equation}
  Combining $(\ast)$ and $(\dagger)$ one gets
  $\Ext_A^{1}(G,\Co[m]{W})=0$. With $g = \RGorpd{\sfU}{M}$ and $m=g+1$
  and $G = \Co[g]{W}$ one now has
  $\Ext_A^{1}(\Co[g]{W},\Co[g+1]{W})=0$. This means that the exact
  sequence $0 \to \Co[g+1]{W} \to W_g \to \Co[g]{W} \to 0$ splits,
  whence $\Co[g]{W}$ is in $\capUV$. Thus one has $\Upd{M} \le g$, and
  the opposite inequality holds by \lemref{dimineq}.
\end{prf*}

\begin{prp}
  \label{prp:ses}
  Let $0 \to M' \to M \to M'' \to 0$ be an exact sequence of
  $A$-modules. With $g' = \RGorpd{\sfU}{M'}$, $g = \RGorpd{\sfU}{M}$,
  and $g'' = \RGorpd{\sfU}{M''}$ there are inequalities,
  \begin{equation*} g' \le \max\{g, g'' - 1\} \,, \quad g \le
    \max\{g', g''\}\,, \quad\text{and}\quad g'' \le \max\{g'+ 1,
    g\}\,.
  \end{equation*}
  In particular, if two of the modules have finite $\catGR$-projective
  dimension, then so has the third.
\end{prp}

\begin{prf*}
  It follows from \prpcite[3.8 and Thm.\ 4.8]{CDEHLT-23} that if two
  of the modules have finite $\catGR$-projective dimension, then so
  does the third. The inequalities are now immediate from
  \rmkcite[4.9]{CDEHLT-23}.
\end{prf*}

\begin{dfn}
  Imitating the Gorenstein global dimension, set
  \begin{equation*}
    \RGorgldim{\sfU} = \sup\setof{\RGorpd{\sfU}{M}}
    { M \text{ is an $A$-module}} \:.
  \end{equation*}
\end{dfn}

The next result shows that the global dimension defined above can be
computed simpler in terms of $\catGR$-projective dimensions of modules
in the subcategory $\sfV$, see also \rmkref{simple}.

\begin{prp}
  \label{prp:gldim}
  The following conditions are equivalent.
  \begin{eqc}
  \item $\RGorgldim{\sfU} < \infty$.
  \item Every $A$-module has finite $\catGR$-projective dimension.
  \item Every $A$-module in $\sfV$ has finite $\catGR$-projective
    dimension.
  \end{eqc}
  Moreover, there is an equality:
  \begin{equation*}
    \RGorgldim{\sfU} = \sup\setof{\RGorpd{\sfU}{M}}
    { M \text{ is an $A$-module in $\sfV$}} \:.
  \end{equation*}
\end{prp}

\begin{prf*}
  The implications
  \eqclbl{i}$\!\implies\!$\eqclbl{ii}$\!\implies\!$\eqclbl{iii} are
  trivial.

  \proofofimp{ii}{i} Assume that every $A$-module has finite
  $\catGR$-projective dimension.  If one has
  $\RGorgldim{\sfU}=\infty$, then for every $n\in\mathbb{N}$ there is
  an $A$-module $M_n$ with $n \leq \RGorpd{\sfU}{M_n} < \infty$.  The
  $A$-module $M = \coprod_{n\in \mathbb{N}}M_n$ has finite
  $\catGR$-projective dimension, say, $k$. Now
  \rmkcite[4.9]{CDEHLT-23} yields
  \begin{equation*}
    \prod_{n\in\mathbb{N}}\Ext_A^{k+1}(M_n,W) \cong \Ext_A^{k+1}(M, W) = 0
  \end{equation*}
  for every module $W$ in $\capUV$. Now another application of
  \rmkcite[4.9]{CDEHLT-23} yields $\RGorpd{\sfU}{M_n} \le k$ for every
  $n\in \NN$, which is absurd.

  \proofofimp{iii}{ii} Assume that every module in $\sfV$ has finite
  $\catGR$-projective dimension.  For every $A$-module $M$, there is
  an exact sequence
  \begin{equation*}
    0 \lra M \lra V \lra U \lra 0
  \end{equation*}
  of $A$-modules with $V$ in $\sfV$ and $U$ in $\sfU$. As
  $\RGorpd{\sfU}{U} \le 0$ holds by \lemref{dimineq} it is immediate
  from \prpref{ses} that $\RGorpd{\sfU}{M} \leq \RGorpd{\sfU}{V}$
  holds. This also proves the asserted equality.
\end{prf*}

To parse the next result, recall the notion of right periodicity from
\dfnref{periodic}. For the absolute cotorsion pair $(\Prj,\Mod)$ the
fact that \eqclbl{i} below implies \eqclbl{ii} is known from
\thmcite[VII.2.2]{ABlIRt07}, but for $n > 0$ the converse appears to
be new.

\begin{thm}
  \label{thm:8-UV}
  Assume that $\sfUV$ is \wrightperiodic\ and let $n \ge 0$ be an
  integer.  The following conditions are equivalent.
  \begin{eqc}
  \item $\RGorgldim{\sfU} \leq n$.
  \item ${\sfU}_n^\perp = \GI$.
  \end{eqc}
  Moreover, when these conditions are satisfied, there is an equality:
  \begin{equation*}
    \RGorgldim{\sfU} = \relFPD{\sfU} \:.
  \end{equation*}
\end{thm}

\begin{prf*}
  Assume that $(i)$ holds. We first argue that one has
  $\relFPD{\sfU}\leq n$. For an $A$-module $M$ with $\Upd M$ finite,
  \lemref{dimineq} yields
  \begin{equation*}
    \Upd M = \RGorpd{\sfU}{M} \leq n \:.
  \end{equation*}
  Next we argue that $\relspli{\sfU}\leq n$ holds. For an injective
  $A$-module $I$, one has per \thmref{SymU} the equality
  $\Upd I=\RGorpd{\sfU}{I} \leq n$. Now ${\sfU}_n^\perp = \GI$ holds
  by \prpref{4-UV} as $\sfUV$ is \wrightperiodic.

  Now assume that \eqclbl{ii} holds. \prpref{4-UV} yields
  \mbox{$\relspli{\sfU} = \relFPD{\sfU} \leq n$}.  By \prpref{gldim}
  it suffices to show that every module in $\sfV$ has
  $\mathsf{RGor_U}(A)$-projective dimension at most $n$. Let $V$ be a
  module in $\sfV$. As $\sfV$ is coresolving, there is an exact sequence
  $0 \to V \to I \to V' \to 0$ with $I$ injective and $V'$ in
  $\sfV$. Taking successive special $\sfU$-precovers one gets exact
  sequences,
  \begin{gather*}
    \tag{$\sharp$} \cdots \xra{\partial} W_n \lra W_{n-1} \lra \cdots
    \lra W_0 \lra V \lra 0 \intertext{and}
    \tag{$\sharp'$} \cdots \xra{\partial'} W_n' \lra W_{n-1}' \lra
    \cdots \lra W_0' \lra V' \lra 0 \:,
  \end{gather*}
  and the Horseshoe Lemma \lemcite[8.2.1]{RHA} yields the exact
  sequence,
  \begin{equation*}
    \cdots \xra{\eth} W_n \oplus W_n' \lra W_{n-1} \oplus W_{n-1} '
    \lra W_{0}\oplus W'_0 \lra I \lra 0 \:.
  \end{equation*}
  With $G = \Coker(\partial)$, $F_{-1} = \Coker(\eth)$, and
  $G' = \Coker(\partial')$ one gets an exact sequence
  $0 \to G \to F_{-1} \to G' \to 0$. By construction $G$ and $G'$
  belong to $\sfV$ and, therefore, so does $F_{-1}$. It follows from
the already established inequality $\relspli{\sfU} \leq n$ and \prpref{Upd} that
$F_{-1}$ belongs to $\sfU$. Repeating this process, one establishes an
exact sequence
  \begin{equation*}
    \tag{$\ast$}
    0 \lra G \lra F_{-1} \lra F_{-2} \lra \cdots
  \end{equation*}
  where the modules $F_i$ are in $\capUV$. Taking successive special
  $\sfU$-precovers yields another exact sequence
  \begin{equation*}
    \tag{$\ast\ast$}
    \cdots \lra F_1 \lra F_0 \lra G \lra 0
  \end{equation*}
  of $A$-modules with each $F_i$ in $\capUV$. Splicing the sequences
  $(\ast)$ and $(\ast\ast)$, one gets an acyclic complex $F$ of
  modules in $\capUV$, where each cycle module $\Cy[i]{F}$ belongs to
  $\sfV$ by construction, and $G = \Co[0]{F}$. As every module in
  $\capUV$ by \prpref{7-UV} has finite injective dimension, at most
  $n$, the complex $\Hom(F,W)$ is acyclic for every module $W$ in
  $\capUV$. Thus, $F$ is a right $\sfU$-totally acyclic complex and
  $G$ is in $\mathsf{RGor_U}(A)$, so from $(\sharp)$ one gets
  $\RGorpd{\sfU} V \le n$.
\end{prf*}

Applied to the absolute cotorsion pair $(\Prj,\Mod)$, the equivalence
of \eqclbl{i} and \eqclbl{iii} in the next theorem recovers the
characterization of left Gorenstein rings from
\corcite[VII.2.6]{ABlIRt07}.

\begin{thm}
  \label{thm:9-1-UV}
  Assume that $\sfUV$ is \wrightperiodic\ and let $n \ge 0$ be an
  integer. The following conditions are equivalent.
  \begin{eqc}
  \item $\RGorgldim{\sfU} \leq n$.
  \item $\sup\setof{\Gid M}{M \in \sfV} \leq n$.
  \item $\sfU_n \cap \sfV = \Inj_n \cap \sfV$.
  \item $\relspli{\sfU} = \sup\setof{\id M}{M \in \capUV} \leq n$.
  \end{eqc}
\end{thm}

\begin{prf*}
  \proofofimp{i}{ii} It follows from \thmref{8-UV} and \fctref{Un}
  that $({\sfU}_n,\GI)$ is a complete cotorsion pair. Let $V$ be a
  module in $\sfV$. There is an exact sequence,
  \begin{equation*}
    0 \lra G \lra U \lra V \lra 0 \:,
\end{equation*}
of $A$-modules with $G\in \GI$ and $U \in {\sfU}_n$.  As $\sfUV$ is
\wrightperiodic, every Gorenstein injective $A$-module
belongs to $\sfV$, so the module $U$ is in $\sfV$ and
$\id U \leq \max\{\relspli{\sfU}, \relFPD{\sfU} \}$ holds by
\prpref{7-UV}. From \prpref{4-UV} and \thmref{8-UV} one now gets
$\relspli{\sfU} = \relFPD{\sfU} \leq n$. Consequently, one has
$\Gid V \leq n$ by Holm \thmcite[2.22 and 2.25]{HHl04a}.

  \proofofimp{ii}{iii} Let $M\in \sfU_n\cap \sfV$. By assumption,
  $\Gid M\leq n$ holds, and so there is an exact sequence
  $0 \to G \to E \to M \to 0$ with $\id E\leq n$ and $G$ Gorenstein
  injective; see Christensen, Frankild, and Holm
  \lemcite[2.18]{CFH-06}. The sequence splits as $G$ belongs to
  $\sfU_n^\perp$, see \prpref{periodic-n}, thus $\id M\leq n$. For the
  opposite containment, suppose that $M$ belongs to $\sfV$ and
  $\id M\leq n$ holds. Let $V\in \sfV$; as $\Gid V\leq n$ holds,
  \thmcite[2.22]{HHl04a} yields $\Ext_A^{n+1}(M,V)=0$, that is,
  $\Upd M\leq n$ holds by \prpref{Upd}.

  \proofofimp{iii}{iv} Set $s = \relspli{\sfU}$ and
  $t = \sup\setof{\id M}{M \in \capUV}$. As every injective module is
  in $\sfV$, the assumption immediately yields $s \le n$. Similarly,
  $t$ is finite, at most $n$. Now there is an injective $A$-module $I$
  with $\Upd I = s$, so taking successive special $\sfU$-precovers one
  gets an exact sequence,
  \begin{equation*}
    0 \lra U_s \lra U_{s-1} \lra \cdots \lra U_0 \lra I \lra 0 \:,
  \end{equation*}
  of $A$-modules with each module $U_i$ in \mbox{$\capUV$}. One has
  $\Ext_A^{s}(I,U_s)\neq 0$ and, therefore, $\id U_s \geq s$, so one
  has $t \geq s$.  On the other hand, there is a module $M$ in
  $\capUV$ with $\id M =t$.  So there exists an exact sequence
  \begin{equation*}
    0 \lra M \lra I_{0} \lra I_{-1} \lra \cdots \lra I_{-t} \lra 0
  \end{equation*}
  of $A$-modules with each $I_i$ injective. One has
  $\Ext_A^{t}(I_{-t}, M)\neq 0$ and, therefore, $\Upd{I_{-t}} \geq t$,
  so one has $s \geq t$.

  \proofofimp{iv}{i} Let $M$ be an $A$-module of finite
  $\sfU$-projective dimension.  If $M$ is in $\sfU$, then
  $\Upd{M} \le n$ holds trivially. Now assume that $d = \Upd{M}$ is
  positive and consider the exact sequence,
  \begin{equation*}
    0 \lra U_d \lra U_{d-1} \lra \cdots \lra U_0 \lra M \lra 0 \:,
  \end{equation*}
  obtained by taking special $\sfU$-precovers. For $i > 0$ the module
  $U_i$ is in $\capUV$. One has $\Ext_A^{d}(M,U_d) \neq 0$ and since,
  by assumption, $\id{U_d} \le n$ holds, one gets $d \le n$.  Thus
  $\relFPD{\sfU} \le n$ holds, and now \prpref{4-UV} and \thmref{8-UV}
  yield the desired inequality.
\end{prf*}

Applied to the absolute cotorsion pair $(\Prj,\Mod)$ the next
corollary recovers Bennis and Mahdou's result \thmcite[1.1]{DBnNMh10}
that the Gorenstein global dimension can be computed based on the
Gorenstein injective dimension.

\begin{cor}
  \label{cor:9-1-UV}
  Assume that $\sfUV$ is \wrightperiodic. The $\catGR$-global dimension
  is finite if and only if every module in $\sfV$ has finite
  Gorenstein injective dimension, and there is an equality:
  \begin{equation*}
    \RGorgldim{\sfU} = \sup\setof{\Gid M}{M \in \sfV} \:.
  \end{equation*}
\end{cor}

\begin{prf*}
  The equality is immediate from the equivalence of conditions
  \eqclbl{i} and \eqclbl{ii} in \thmref{9-1-UV}. The first assertion
  also follows from \thmref{9-1-UV} as $\sfV$ is closed under
  products.
\end{prf*}

For the remainder of the section, we shift focus to the class $\sfV$,
recalling the notion of left $\sfV$-Gorenstein modules from
\dfncite[2.1]{CET-20} and the associated homological dimension
$\LGorid{\sfV}$ from \dfncite[4.5]{CDEHLT-23}.  For
$(\sfU,\sfV) = (\Mod,\Inj)$ the $\catGL$-injective dimension is the
standard Gorenstein injective dimension, and we use the standard
notation $\mathrm{Gid}_A$, as already done in \thmref{9-1-UV} and
\corref{9-1-UV}.

\begin{lem}
  \label{lem:dimineq-bis}
  Let $M$ be an $A$-module. There is an inequality,
  \begin{equation*}
    \LGorid{\sfV}{M} \leq \Vid M \:,
  \end{equation*}
  and equality holds if\, $\Vid M$ is finite.
\end{lem}

\begin{prf*}
  Dual to the proof of \lemref{dimineq}.
\end{prf*}

\begin{thm}
  \label{thm:SymV}
  Let $M$ be an $A$-module. If\, $\pd M$ is finite, then one has
  \begin{equation*}
    \LGorid{\sfV}{M} = \Vid M \:.
  \end{equation*}
\end{thm}

\begin{prf*}
  Dual to the proof of \thmref{SymU}.
\end{prf*}

\begin{prp}
  \label{prp:ses-bis}
  Let $0 \to M' \to M \to M'' \to 0$ be an exact sequence of
  $A$-modules. With $g' = \LGorid{\sfV}{M'}$, $g = \LGorid{\sfV}{M}$,
  and $g'' = \LGorid{\sfV}{M''}$ there are inequalities,
  \begin{equation*}
    g' \le \max\{g, g'' + 1\} \,, \quad g \le \max\{g',g''\}\,,
    \quad\text{and}\quad g'' \le \max\{g'- 1, g\}\,.
  \end{equation*}
  In particular, if two of the modules have finite $\catGL$-injective
  dimension, then so has the third.
\end{prp}

\begin{prf*}
  Dual to the proof of \prpref{ses}.
\end{prf*}

\begin{dfn}
  Set
  \begin{equation*}
    \LGorgldim{\sfV} = \sup\setof{\LGorid{\sfV}{M}}
    { M \text{ is an $A$-module}} \:.
  \end{equation*}
  For $\sfUV = (\Mod,\Inj)$ this is by \thmcite[1.1]{DBnNMh10} the
  invariant $\Ggldim$.
\end{dfn}

\begin{prp}
  \label{prp:gldim-bis}
  The following conditions are equivalent.
  \begin{eqc}
  \item $\LGorgldim{\sfV} < \infty$.
  \item Every $A$-module has finite $\catGL$-injective dimension.
  \item Every $A$-module in $\sfU$ has finite $\catGL$-injective
    dimension.
  \end{eqc}
  Moreover, there is an equality:
  \begin{equation*}
    \LGorgldim{\sfV} = \sup\setof{\LGorid{\sfV}{M}}{ M \text{ is an $A$-module in $\sfU$}} \:.
  \end{equation*}
\end{prp}

\begin{prf*}
  Dual to the proof of \prpref{gldim}.
\end{prf*}

\begin{thm}
  \label{thm:8-UV-dual}
  Assume that $\sfUV$ is \wleftperiodic\ and let $n \ge 0$ be an
  integer.  The following conditions are equivalent.
  \begin{eqc}
  \item $\LGorgldim{\sfV}\leq n$.
  \item ${{}^\perp\sfV}_n = \GP$.
  \end{eqc}
  Moreover, when these conditions are satisfied, there is an equality:
  \begin{equation*}
    \LGorgldim{\sfV} = \relFID{\sfV} \:.
  \end{equation*}
\end{thm}

\begin{prf*}
  Dual to the proof of \thmref{8-UV}.
\end{prf*}

\begin{thm}
  \label{thm:9-1-UV-bis}
  Assume that $\sfUV$ is \wleftperiodic\ and let $n \ge 0$ be an
  integer. The following conditions are equivalent.
  \begin{eqc}
  \item $\LGorgldim{\sfV} \leq n$.
  \item $\sup\setof{\Gpd M}{M \in \sfU} \leq n$.
  \item $\sfU \cap \sfV_n = \sfU \cap \Prj_n$.
  \item $\relsilp{\sfV} = \sup\setof{\pd M}{M \in \capUV} \leq n$.
  \end{eqc}
\end{thm}

\begin{prf*}
  Dual to the proof of \thmref{9-1-UV}.
\end{prf*}

\begin{cor}
  \label{cor:9-1-UV-bis}
  Assume that $\sfUV$ is \wleftperiodic. The $\catGL$-global dimension
  is finite if and only if every module in $\sfU$ has finite
  Gorenstein projective dimension, and there is an equality:
  \begin{equation*}
    \LGorgldim{\sfV} = \sup\setof{\Gpd M}{M \in \sfU} \:.
  \end{equation*}
\end{cor}

\begin{prf*}
  Dual to the proof of \corref{9-1-UV}.
\end{prf*}

We denote by $\Flat$ the class of flat $A$-modules and by $\Cot$ the
class of cotorsion $A$-modules; they form a cotorsion pair.

\begin{rmk}
  \label{rmk:simple}
  The equalities in \prpref{gldim} and \corref{9-1-UV}, and dually in
  \prpref{gldim-bis} and \corref{9-1-UV-bis}, are the counterparts in
  Gorenstein homological algebra of the following equalities:
  \begin{align*}
    \sfU\operatorname{-gldim}(A) = \sup\setof{\Upd V}{V \in \sfV}
    = \sup\setof{\id V}{V \in \sfV} \:, \\
    \sfV\operatorname{-gldim}(A) =\sup\setof{\Vid U}{U \in \sfU}
    = \sup\setof{\pd U}{U \in \sfU} \:.
  \end{align*}
  Here, of course, $\sfU\operatorname{-gldim}(A)$ and
  $\sfV\operatorname{-gldim}(A)$ are defined as the suprema of $\Upd$
  and $\Vid$ over all $A$-modules. The displayed equalities can be
  shown as in Mao and Ding \thmcite[19.2.5(1)]{LMaNDn06}, where they
  consider the equalities in the second line for the cotorsion pair
  $(\sfU,\sfV)=(\Flat,\Cot)$.
\end{rmk}


\section{Gorenstein flat-cotorsion global dimension}
\label{sec:gfcgldim}

\noindent
In this section we apply the theory developed in the previous two
sections to the cotorsion pair $(\Flat,\Cot)$ in order to study the
global invariant based on the Gorenstein flat-cotorsion dimension from
\cite{CELTWY-21}. We relate it to the Gorenstein global dimension and
obtain new characterizations of Iwanaga--Gorenstein rings and
Ding--Chen rings. The results in this section find further
applications in \secref{gwgldim}.

The cotorsion pair $(\Flat,\Cot)$ is hereditary and generated by a
set; see Bican, El Bashir, and Enochs \prpcite[2]{BEE-01}. The
relevant special case of \fctref{Un} was noted earlier by H\"ugel,
Herbera, and Trlifaj \seccite[1.C]{HHT-06}. The cotorsion pair
$(\Flat,\Cot)$ is \rightperiodic\ by \thmcite[1.3]{BCE-20}.

In this setting, $\sfU = \Flat$, the invariants $\Upd{\!}$,
$\relFPD{\sfU}$, and $\relspli{\sfU}$ are the flat dimension
$\fd{\!}$, the finitistic flat dimension $\FFD$, and $\sfli$. Modules
in $\catGR(A)$ are called Gorenstein flat-cotorsion, see
\rmkcite[4.8]{CELTWY-21}. The $\catGR$-projective dimension is known
as the Gorenstein flat-cotorsion dimension and written
$\Gfcd{\!}$. For the corresponding global dimension we naturally
introduce the symbol $\Gfcgldim[A]$.

\begin{fct}
  \label{fct:gwgldim}
  It follows from \thmcite[5.7]{CELTWY-21} that there is an
  inequality,
  \begin{equation*}
    \Gfcgldim \le \Gwgldim \:,
  \end{equation*}
  and equality holds if $\Gwgldim$ is finite.
\end{fct}

For a large class of rings the Gorenstein flat-cotorsion global
dimension agrees with the Gorenstein weak global dimension, see
\secref{gwgldim}. As such, the next theorem can be considered an
improvement of \thmcite[5.3]{IEm12}, see also
\thmcite[2.4]{CET-21a}. Indeed, $\Gwgldim$ is finite if and only if
$\sfli$ and $\sfli[\Aop]$ are finite, but the next theorem yields
\begin{equation*}
  \Gfcgldim = \max\{\sfli, \FFD\} \:,
\end{equation*}
and one always has $\FFD \le \sfli[\Aop]$ by \prpcite[2.4]{IEmOTl11}.

\begin{thm}
  \label{thm:8}
  Let $n \ge 0$ be an integer. The following conditions are
  equivalent.
  \begin{eqc}
  \item $\Gfcgldim[A]\leq n$.
  \item $\max\{\sfli, \FFD\} \leq n$.
  \item $\sup\setof{\Gid{M}}{M \text{ is cotorsion}} \le n$.
  \item $\Flat_n \cap \Cot = \Inj_n \cap \Cot$.
  \item
    $\sfli = \sup\setof{\id{M}}{M \text{ is flat-cotorsion}} \le n$.
  \item $\Flat_n^\perp = \GI$.
  \end{eqc}
  Moreover, when these conditions are satisfied there are equalities:
  \begin{equation*}
    \Gfcgldim = \FFD = \sfli\:.
  \end{equation*}
\end{thm}

\begin{prf*}
  Apply \prpref{4-UV} and \thmref[Theorems~]{8-UV} and
  \thmref[]{9-1-UV} to $(\Flat, \Cot)$.
\end{prf*}

The next corollary provides the converse to \prpcite[4.2]{CET-23} that
was anticipated in \cite[Quest.~4.11]{CET-23}.  Recall that an
$A$-module $M$ is said to be \emph{strongly cotorsion} if
$\Ext_A^1(L,M)=0$ for every $A$-module $L$ of finite flat dimension.

\begin{cor}
  \label{cor:question}
  The following conditions are equivalent.
  \begin{eqc}
  \item $\Gfcgldim < \infty$.
  \item $\FFD < \infty$ and every strongly cotorsion $A$-module is
    Gorenstein injective.
  \end{eqc}
\end{cor}

\begin{prf*}
  With $n = \FFD$ the class of strongly cotorsion modules is
  $\Flat_n^\perp$ and the assertion follows from \thmref{8}.
\end{prf*}

As discussed in \seccite[4]{CET-23} there has been interest in the
question of when Gorenstein injective modules coincide with strongly
cotorsion modules.

\begin{rmk}
  Assume that $\FFD=n$ is finite and let $\SC$ be the class of
  strongly cotorsion $A$-modules. The pair $(\Flat_n,\SC)$ is per
  \fctref{Un}, or \seccite[1.C]{HHT-06}, a hereditary cotorsion pair
  generated by a set, and it is \rightperiodic\ by
  \prpref{periodic-n}. Applied to this cotorsion pair, \corref{9-1-UV}
  yields
  \begin{equation*}
    \RGorgldim{\Flat_n}= \sup\setof{\Gid M}
    {M \text{ is a strongly cotorsion $A$-module}} \:.
  \end{equation*}
  This measures how far the strongly cotorsion $A$-modules are from
  being Gorenstein injective: It is zero precisely when the strongly
  cotorsion and Gorenstein injective modules coincide, or
  alternatively when $\Gfcgldim$ is finite; see \corref{question}.
\end{rmk}

The next corollary is subsumed by \thmref{gfcgf-gldim}, but we record
it here for use in \corref{IG}, \corref{DingChen}, and \thmref{10}.

\begin{cor}
  \label{cor:gwgldim}
  If $A$ is left or right coherent, then there is an equality:
  \begin{equation*}
    \Gfcgldim = \Gwgldim \:.
  \end{equation*}
\end{cor}

\begin{prf*}
  By \fctref{gwgldim} one has the inequality $\Gfcgldim \leq \Gwgldim$
  with equality if $\Gwgldim$ is finite.  It follows from
  \corref{question} and \thmcite[4.9]{CET-23} that the two global
  dimensions are simultaneously finite over a left or right coherent
  ring, whence the displayed equality holds.
\end{prf*}

\begin{prp}
  \label{prp:9-1}
  The following conditions are equivalent.
  \begin{eqc}
  \item $\Gfcgldim < \infty$.
  \item Every cotorsion $A$-module has finite Gorenstein
    flat-cotorsion dimension.
  \item Every cotorsion $A$-module has finite Gorenstein injective
    dimension.
  \end{eqc}
  Further, there are equalities,
  \begin{align*}
    \Gfcgldim & = \sup\setof{\Gfcd M}{M \text{ is a cotorsion $A$-module}} \\
              & = \sup\setof{\Gid M}{M \text{ is a cotorsion $A$-module}} \:.
  \end{align*}
\end{prp}

\begin{prf*}
  The equivalence of \eqclbl{i}--\eqclbl{iii} as well as the
  equalities follow from \prpref{gldim} and \corref{9-1-UV} applied to
  the cotorsion pair $(\Flat, \Cot)$.
\end{prf*}

The equalities in \prpref{9-1} are the correspondents in Gorenstein
homological algebra to the fact that the weak global dimension is the
supremum of flat dimensions of cotorsion modules and also the supremum
of injective dimensions of cotorsion modules, see \rmkref{simple}.

\begin{cor}
  \label{cor:IG}
  If $A$ is noetherian, then it is Iwanaga--Gorenstein if (and only
  if) every cotorsion $A$-module has finite Gorenstein injective
  dimension.
\end{cor}

\begin{prf*}
  Combine \corref{gwgldim} and \prpref{9-1} with
  \rmkcite[3.11]{CET-21a}.
\end{prf*}

Recall that $A$ is called \emph{Ding--Chen} if it is coherent with
finite self-fp-injective dimension on both sides. Symmetry of the
Gorenstein weak global dimension has the following immediate
consequence:

\begin{prp}
  \label{prp:DingChen}
  If $A$ is coherent, then it is Ding--Chen if and only if $\Gwgldim$
  is finite.
\end{prp}

\begin{prf*}
  As $A$ is coherent, it follows from work of Ding and Chen
  \thmcite[7]{NDnJCh96} combined with \corcite[2.5]{CET-21a} that $A$
  is Ding--Chen if and only if the Gorenstein weak global dimension is
  finite.
\end{prf*}

One can get an even broader statement than \corref{IG}:

\begin{cor}
  \label{cor:DingChen}
  If $A$ is coherent, then the following conditions are equivalent.
  \begin{eqc}
  \item $A$ is Ding--Chen.
  \item Every cotorsion $A$-module has finite Gorenstein injective
    dimension.
  \item Every $A$-module has finite Gorenstein flat dimension.
  \end{eqc}
\end{cor}

\begin{prf*}
  In view of \thmcite[2.4]{CET-21a} one gets the equivalence of the
  three conditions by combining \corref{gwgldim} with
  \prpref[Propositions~]{9-1} and \prpref[]{DingChen}.
\end{prf*}

To parse the next result recall the invariant
\begin{equation*}
  \splf = \sup\setof{\pd{F}}{F \text{ is a flat $A$-module}} \:.
\end{equation*}
If $A$ is right coherent, the next result reduces per \corref{gwgldim}
to \corcite[3.5]{CET-21a}.

\begin{cor}
  \label{cor:1.3}
  There are inequalities:
  \begin{align*}
    \Gfcgldim \leq \Ggldim \leq \Gfcgldim + \splf \:.
  \end{align*}
  Moreover, the following conditions are equivalent.
  \begin{eqc}
  \item $\Gfcgldim$ and $\splf$ are finite.
  \item $\Ggldim$ is finite.
  \end{eqc}
\end{cor}

\begin{prf*}
  The left-hand inequality holds by \thmcite[3.3]{CET-21a}; it is also
  a consequence of \thmref[Theorems~]{Ggldim-new} and \thmref[]{8}. To
  prove the right-hand inequality one can assume that $g = \Gfcgldim$
  and $s = \splf$ are finite. Let $M$ be an $A$-module and $M \qra I$
  an injective resolution. By dimension shifting or \rmkref{simple}
  the module $\Cy[-s]{I}$ is cotorsion, so by \prpref{9-1} one has
  $\Gid{M} \le \Gid{\Cy[-s]{I}} + s \le g + s$.

  It is clear from the second inequality that \eqclbl{i} implies
  \eqclbl{ii}. For the converse it follows from the first inequality
  that $\Gfcgldim$ is finite and one has $\splf \le \FPD = \Ggldim$ by
  \prpcite[6]{CUJ70} and \thmcite[2.28]{HHl04a}.
\end{prf*}


\section{A Gorenstein analogue of a result due to Stenstr\"om}
\label{sec:stenstrom}

\noindent
Recall that an $A$-module $I$ is called \emph{fp-injective} if
$\Ext_A^1(F,I) =0$ holds for every finitely presented $A$-module $F$.
Let $\FpInj$ be the class of fp-injective $A$-modules.  In a pure
homological sense, fp-injective modules are dual to flat modules as
every exact sequence $0 \to I \to M \to N \to 0$ with $I$ fp-injective
is pure.

\begin{dfn}
  \label{dfn:WGor}
  For simplicity we call $\catGL[\FpInj]$-modules \emph{Gorenstein
    fp-injective} $A$-modules. The associated dimensions
  $\LGorid{\FpInj}{}$ and $\LGorgldim{\FpInj}$ are similarly
  simplified to $\Gfpid{}$ and $\Gfpgldim$.
\end{dfn}

It is standard that $({}^\perp\FpInj,\FpInj)$ is a complete cotorsion
pair generated by the set of isomorphism classes of finitely presented
modules. It is hereditary if and only if $A$ is left coherent; see for
example \thmcite[3.2]{BSt70} for the non-trivial ``only if'' part of
the statement.

\begin{fct}
  \label{fct:fpCycles}
  Let $A$ be left coherent.  Every acyclic complex of $A$-modules in
  ${}^\perp\FpInj$ has cycle submodules in ${}^\perp\FpInj$, that is,
  $({}^\perp\FpInj,\FpInj)$ is \leftperiodic; see \v{S}aroch and
  \v{S}tov{\'{\i}}{\v{c}}ek \exacite[4.3]{JSrJSt20}.
\end{fct}

\begin{rmk}
  An $A$-module is Ding injective if it is a cycle module in an
  acyclic complex of injective modules that remains exact upon
  application of $\Hom(I,-)$ for every fp-injective $A$-module $I$;
  these modules were called Gorenstein fp-injective by Mao and Ding
  \cite{LMaNDn08}. If $A$ is left noetherian, then Ding injective
  $A$-modules and Gorenstein fp-injective $A$-modules are simply
  Gorenstein injective $A$-modules. In general, Ding injective modules
  and Gorenstein fp-injective modules can differ: 
\end{rmk}

\begin{exa}
  Let $A$ be von Neumann regular and not noetherian. As the cotorsion
  pair $(\Flat,\Cot) = (\Mod,\Inj)$ is \rightperiodic, every cycle
  submodule in an acyclic complex of injective $A$-modules is
  injective. Thus, an $A$-module is Ding injective if and only if it
  is injective. Similarly, $({}^\perp\FpInj,\FpInj) = (\Prj,\Mod)$
  holds, so an $A$-module is Gorenstein fp-injective if and only if it
  is projective, see \dfncite[1.1 and Exa.\ 2.2]{CET-20}.
\end{exa}

For $\sfV={\FpInj}$ we use the notation $\fpid$ for $\Vid{\!}$ and,
similarly, $\sfpilp$ for $\relsilp{\sfV}$ and $\FFPID$ for
$\relFID{\sfV}$.

\begin{lem}
  \label{lem:5.5}
  If $A$ is left coherent, then there are equalities,
  \begin{equation*}
    \sfli[\Aop] = \sfpilp[A]
    \quad\text{and}\quad
    \FFD[\Aop] =\FFPID[A] \:.
  \end{equation*}
\end{lem}

\begin{prf*}
  Since $A$ is left coherent, one has $\sfli[\Aop]=\fpid A=\sfpilp[A]$
  by Ding and Chen \thmcite[3.8]{NQDJLC93}.  On the other hand, for
  every $\Aop$-module $M$ and $A$-module $N$, the equalities
  \begin{equation*}
    \fd[\Aop]M = \fpid \Hom[\ZZ](M,\QQ/\ZZ) \quad\text{and}\quad
    \fpid N = \fd[\Aop]\Hom[\ZZ](N,\QQ/\ZZ)
  \end{equation*}
  hold by Fieldhouse \thmcite[2.1 and 2.2]{DJF72}, whence
  $\FFD[\Aop]=\FFPID[A]$.
\end{prf*}

The next result is an analogue in Gorenstein homological algebra to a
result of Stenstr\"om \thmcite[3.3]{BSt70}. The equality of the first
and last quantities can alternatively be deduced from results of
Bennis \thmcite[3.3]{DBn11} and Li, Yan, and Zhang
\thmcite[2.10]{LYZ-23}, but the remaining equalities are new.

\begin{thm}
  \label{thm:10}
  If $A$ is left coherent, then there are equalities,
  \begin{align*}
    \Gwgldim & = \Gfpgldim \\
             & = \sup\setof{\Gfpid M}{M~\text{is in } {}^\perp\FpInj } \\
             & = \sup\setof{\Gpd M}{M~\text{is a finitely
               presented $A$-module}}                \:.
  \end{align*}
\end{thm}

\begin{prf*}
  In view of the equalities in \lemref{5.5} it follows from
  \thmref{8}, \prpref{4-UV-bis}, and \thmref{8-UV-dual} that the
  quantities $\Gfcgldim[\Aop]$ and $\Gfpgldim$ are simultaneously
  finite whence one has
  \begin{equation*}
    \Gwgldim = \Gfcgldim[\Aop] = \FFD[\Aop] =  \FFPID = \Gfpgldim
  \end{equation*}
  by the same references combined with \corcite[2.5]{CET-21a} and
  \corref{gwgldim}. This proves the first of the asserted equalities,
  and the second holds by \prpref{gldim-bis}.
  
  To prove the third equality, notice first that \corref{9-1-UV-bis}
  yields
  \begin{equation*}
    \Gfpgldim = \sup\setof{\Gpd M}{M~\text{is in } {}^\perp\FpInj} \:,
  \end{equation*}
  and since all finitely presented $A$-modules are in
  ${}^\perp\FpInj$, the right-hand quantity is at least
  $\sup\setof{\Gpd M}{M~\text{is a finitely presented $A$-module}}$.
  Assume now that this last quantity is at most $n$ and let $\sfS$ be
  a set of representatives of isomorphism classes of finitely
  presented $A$-modules. Let $\sfL$ be the class of $A$-modules of
  Gorenstein projective dimension at most $n$.  By assumption,
  $\sfS\subseteq \sfL$ and, therefore,
  $\Filt(\sfS)\subseteq \Filt(\sfL)$. Here, for a class $\sfC$ of
  $A$-modules, $\Filt(\sfC)$ is the class of direct transfinite
  extensions of modules in $\sfC$. By Enochs, Iacob, and Jenda
  \thmcite[3.4]{EIJ-07} one has $\Filt(\sfL)\subseteq \sfL$, so all
  modules in $\Filt(\sfS)$ have Gorenstein projective dimension at
  most $n$. Since $\sfS$ is a set that contains $A$,
  \corcite[6.14]{RGbJTr12-1} yields that
  ${}^\perp(\sfS^\perp)=\add(\Filt(\sfS))$. It is standard that $\sfL$
  is closed under direct summands. Since $({}^\perp\FpInj,\FpInj)$ is
  generated by the set $\sfS$, one has
  ${}^\perp\FpInj={}^\perp(\sfS^\perp)$, so all modules in
  ${}^\perp\FpInj$ have Gorenstein projective dimension at most $n$.
\end{prf*}

The next corollary is an equivalent in Gorenstein homological algebra
to \thmcite[4]{MAs55}. As noticed in \corcite[2.7]{LYZ-23} one could
alternatively derive it by combining a result of Bouchiba
\thmcite[7]{SBc15a} with symmetry of the Gorenstein weak global
dimension \corcite[2.5]{CET-21a}. As also noticed in \cite{LYZ-23}
this symmetry could alternatively be obtained by combining
\thmcite[6]{SBc15a} with \corcite[4.12]{JSrJSt20}.

\begin{cor}
  \label{cor:Bouchiba}
  If $A$ is left noetherian, then the following equality holds:
  \begin{equation*}
    \Gwgldim = \Ggldim \:.
  \end{equation*}
\end{cor}

\begin{prf*}
  As $A$ is left noetherian, fp-injective modules are simply injective
  modules and one has $\Gfpid{M} = \Gid{M}$ for every $A$-module $M$.
  Now invoke the first equality in \thmref{10}.
\end{prf*}

\begin{cor}
  If $A$ is coherent, then it is Ding--Chen if and only if every
  $A$-module has finite Gorenstein fp-injective dimension.
\end{cor}

\begin{prf*}
  Combine \prpref{DingChen} with \prpref{gldim-bis} and \thmref{10}.
\end{prf*}

We observe that the left periodicity of $({}^\perp\FpInj,\FpInj)$
provides for a nicer description of Gorenstein fp-injective modules.

\begin{prp}
  \label{prp:fpinj}
  Let $A$ be left coherent. An $A$-module $G$ is Gorenstein
  fp-injective if and only if there exists an acyclic complex $T$ of
  modules from the class \mbox{$\sfW= {}^\perp\FpInj \cap \FpInj$}
  with $G = \Cy[0]{T}$ such that $\Hom(W,T)$ is acyclic for every
  module $W$ in $\sfW$.
\end{prp}

\begin{prf*}
  The assertion follows from the definition of Gorenstein fp-injective
  modules, as the left periodic property of $({}^\perp\FpInj,\FpInj)$
  takes care of condition (3).
\end{prf*}

\begin{rmk}
  Assuming that $A$ is not left coherent,
  $({}^\perp\FpInj,\FpInj)$ is not a hereditary cotorsion pair.  To obtain a cotorsion
  pair that is hereditary, generated by a set, and agrees with
  $({}^\perp\FpInj,\FpInj)$ in the left coherent case, one can replace
  the fp-injective modules by the smaller class of strongly
  fp-injective modules. This yields a cotorsion pair which is
  \wleftperiodic\ by Emmanouil and Kaperonis
  \rmkcite[4.10(ii)]{IEmIKp}, and we can thus apply the results of
  \secref[Sections~]{uv} and \secref[]{relgor} to this cotorsion
  pair. The question of left periodicity of this cotorsion pair has
  been recently considered by Bazzoni, Hrbek, and Positselski
  \cite{BHP}.
\end{rmk}


\section{Applications to the Gorenstein (weak) global dimension}
\label{sec:gwgldim}

\noindent
In settings where the Gorenstein flat-cotorsion global dimension
agrees with the Gorenstein weak global dimension, results from the
previous sections yield new relations among the Gorenstein global
dimensions, see \thmref[]{GwGgl_ineq}--\thmref[]{Osofsky}. It was
established in \corcite[5.8]{CELTWY-21} that the equality
$\Gfcgldim = \Gwgldim$ holds if $A$ is right coherent, and by
\corref{gwgldim} it also holds if $A$ is left coherent. In this
section, we enlarge the class of rings for which the equality holds in
three directions, see \thmref[Theorems~]{gfcgf-gldim}, \thmref[]{58},
and \thmref[]{a0}.

For an $A$-module $M$ we use the abbreviated notation $M^+$ for the
character module $\Hom[\ZZ](M,\QQ/\ZZ)$. The definable closure of $A$
is denoted $\langle A \rangle$.  Recall from Cort\'es Izurdiaga
\cite{MCI16} that a ring $A$ is \emph{right weak coherent} if the
product of flat $A$-modules has finite flat dimension.

\begin{lem}
  \label{lem:n}
  If $A$ is right weak coherent, then there are inequalities:
  \begin{equation*}
    \sup\setof{\fd{I^+}}{I\text{ is an injective $\Aop$-module}}
    \leq \sup\setof{\fd{M}}{M\in \langle A\rangle} < \infty \:.
  \end{equation*}
\end{lem}

\begin{prf*}
  The character module of every injective $\Aop$-module belongs to
  $\langle A \rangle$, see Estrada, Iacob, and P\'erez
  \rmkcite[2.12]{EIP-20}, whence the left-hand inequality holds. To
  see that the right-hand supremum is finite, notice that the
  definable closure of $A$ agrees with the definable closure of the
  class of flat $A$-modules.  By a result of Rothmaler, see Herbera
  \lemcite[2.9]{DHr14}, the definable closure of the class of flat
  $A$-modules can be constructed by first closing under direct
  products, then direct limits and, finally, pure submodules. As $A$
  is right weak coherent, there is per \prpcite[4.1]{MCI16} an integer
  $n$ such that the product of any family of flat $A$-modules has flat
  dimension at most $n$. The class of modules of flat dimension at
  most $n$ is closed both under direct limits and pure submodules, so
  all modules in $\langle A \rangle$ have finite flat dimension, and
  so the right-hand supremum is finite, at most $n$.
\end{prf*}

The next theorem subsumes \thmcite[5.2]{CELTWY-21}.

\begin{thm}
  \label{thm:gfcgf-1}
  Let $M$ be an $A$-module. If $M$ is Gorenstein flat and cotorsion,
  then it is Gorenstein flat-cotorsion. The converse holds if $A$ is
  right weak coherent.
\end{thm}

\begin{prf*}
  Per \thmcite[5.2]{CELTWY-21} an $A$-module that is both Gorenstein
  flat and cotorsion is Gorenstein flat-cotorsion, so it suffices to
  show that the converse holds under the assumption that $A$ is right
  weak coherent.  Let $M$ be a Gorenstein flat-cotorsion $A$-module
  and $T$ a totally acyclic complex of flat-cotorsion $A$-modules with
  $M=\Co[0]{T}$. As $M$ is cotorsion, see e.g.\
  \lemcite[3.2]{CELTWY-21}, it remains to show that it is Gorenstein
  flat. Let $I$ be an injective $\Aop$-module. The $A$-module $I^+$ is
  cotorsion, and since $A$ is right weak coherent, it follows from
  \lemref{n} that $I^+$ has finite flat dimension. By a standard
  dimension shifting argument, the complex $\Hom(T,I^+)$ is
  acyclic. Adjunction now yields
  \begin{equation*}
    \Hom[\ZZ](I\otimes_AT,\QQ/\ZZ)\cong \Hom(T,I^+) \:,
  \end{equation*}
  and so $I\otimes_AT$ is acyclic by faithful injectivity of
  $\QQ/\ZZ$. Thus $T$ is \textbf{F}-totally acyclic, and $M$ is
  Gorenstein flat.
\end{prf*}

\begin{thm}
  \label{thm:gfcgf-dim}
  Let $A$ be right weak coherent. For every $A$-complex $M$ one has:
  \begin{equation*}
    \Gfcd M=\Gfd M \:.
  \end{equation*}
\end{thm}

\begin{prf*}
  The inequality $\Gfcd M\leq \Gfd M$ holds by
  \thmcite[5.7]{CELTWY-21}, and equality holds trivially if $\Gfcd M$
  is infinite. Assume now that $\Gfcd M=n<\infty$ holds. By definition
  there exists a semi-flat-cotorsion replacement $W$ of $M$ with
  $\H[i]{W}=0$ for $i>n$ and $\Co[n]{W}$ a Gorenstein flat-cotorsion
  $A$-module. By \thmref{gfcgf-1}, the module $\Co[n]{W}$ is
  Gorenstein flat, whence $\Gfd M\leq \Gfcd M$ holds by the definition
  of Gorenstein flat dimension.
\end{prf*}

\begin{rmk}
  While right weak coherence is a natural generalization of right
  coherence, notice that the previous two theorems remain true under
  the assumption that the character module of every injective
  $\Aop$-module has finite flat dimension.
\end{rmk}

\begin{thm}
  \label{thm:gfcgf-gldim}
  If $A$ is left or right weak coherent, then the next equality holds:
  \begin{equation*}
    \Gfcgldim = \Gwgldim \:.
  \end{equation*}
\end{thm}

\begin{prf*}
  If $A$ is right weak coherent, then the equality follows immediately
  from \thmref{gfcgf-dim}. Assume now that $A$ is left weak coherent
  and notice that by \fctref{gwgldim} one has the inequality
  $\Gfcgldim \leq \Gwgldim$ with equality if $\Gwgldim$ is finite.  If
  $\Gfcgldim$ is infinite, then the equality trivially holds, so we
  assume that $\Gfcgldim$ is finite, and proceed to show that
  $\Gwgldim$ is finite. It suffices per \corcite[2.5]{CET-21a} to show
  $\Gwgldim[\Aop]$ is finite.

  As $\Aop$ is right weak coherent, it follows from \lemref{n} that
  there is an integer $n$ with
  \begin{equation}
    \tag{$\ast$}
    \sup\setof{\fd[\Aop]{I^+}}{I\text{ is an injective $A$-module}}\leq n
    \quad\text{and}\quad \Gfcgldim\leq n \:.
  \end{equation}
  Let $M$ be an $\Aop$-module. As $M^+$ is a cotorsion $A$-module,
  \prpref{9-1} yields $\Gid{M^+}\leq n$. By \lemcite[2.18]{CFH-06}
  there is an exact sequence of $A$-modules,
  \begin{equation*}
    0 \lra H \lra E \lra M^+ \lra 0 \:,
  \end{equation*}
  with $\id{E}\leq n$ and $H$ is Gorenstein injective. Taking
  character modules yields an embedding $M^{++} \to E^+$ and,
  therefore, an embedding $M \to E^+$ of $\Aop$-modules. Since
  $\id{E}\leq n$ holds one has $\fd[\Aop]{E^+}\leq 2n$ per the first
  inequality in $(\ast)$. Thus every $\Aop$-module can be embedded
  into an $\Aop$-module of flat dimension at most $2n$. This process
  allows one to build a right resolution of $M$ by $\Aop$-modules
  $M_i$ of flat dimension at most $2n$.  Following Cartan and
  Eilenberg \cite[Chap.~XVII, \S1]{CarEil} one can construct a
  projective resolution of this complex in the category of
  $\Aop$-complexes:
  \begin{equation*}
    \xymatrix@=1.5pc{
      & \vdots\ar[d] & \vdots\ar[d] & \vdots\ar[d] &  \\
      0 \ar[r] & P_1 \ar[r]\ar[d] & P_1^{(0)} \ar[r] \ar[d]& P_1^{(-1)} \ar[r]\ar[d] & \ \cdots\\
      0 \ar[r] & P_0 \ar[r]\ar[d] & P_0^{(0)} \ar[r]\ar[d] & P_0^{(-1)} \ar[r]\ar[d] & \ \cdots\\
      0 \ar[r] & M \ar[r] \ar[d] & M_0 \ar[r] \ar[d] & M_{-1} \ar[r] \ar[d] & \ \cdots\\
      & 0 & 0 & 0 &
    }
  \end{equation*}
  This induces an exact sequence
  \begin{equation*}
    0 \lra \Co[2n]{P} \lra \Co[2n]{P^{(0)}} \lra \Co[2n]{P^{(-1)}} \lra \cdots \:.
  \end{equation*}
  As the modules $\Co[2n]{P^{(i)}}$ are flat, this means that a
  syzygy, namely $\Co[2n]{P}$, of $M$ is a cokernel in an acyclic
  complex of flat $\Aop$-modules. The inequality $\Gfcgldim \leq n$
  implies $\sfli\leq n$ by \thmref{8}. Thus every acyclic complex of
  flat $\Aop$-modules is \textbf{F}-totally acyclic, whence $M$ has
  finite Gorenstein flat dimension.
\end{prf*}

\begin{cor}
  \label{thm:gfcgf-sym}
  If $A$ is left or right weak coherent, then the next equality holds:
  \begin{equation*}
    \Gfcgldim[A]=\Gfcgldim[\Aop] \:.
  \end{equation*}
\end{cor}

\begin{prf*}
  As both $A$ and $\Aop$ are left or right weak coherent, the equality
  follows from \thmref{gfcgf-gldim} and \corcite[2.5]{CET-21a}.
\end{prf*}

One also obtains the following extension of \thmcite[4.9]{CET-23}:

\begin{cor}
  If $A$ is left or right weak coherent, then the following conditions
  are equivalent.
  \begin{eqc}
  \item $\Gwgldim < \infty$.
  \item $\FFD < \infty$ and every strongly cotorsion $A$-module is
    Gorenstein injective.
  \end{eqc}
\end{cor}

\begin{prf*}
  This is immediate from \corref{question} in light of
  \thmref{gfcgf-gldim}.
\end{prf*}

Per \prpcite[6]{CUJ70} it follows from next theorem that the equality
$\Gfcgldim = \Gwgldim$ holds if $\FPD$ is finite.

\begin{thm}
  \label{thm:58}
  If $\Cot\operatorname{-id}_{A}A$ is finite, then the following
  equality holds:
  \begin{equation*}
    \Gfcgldim = \Gwgldim \:.
  \end{equation*}
  In particular, the equality holds if $\splf$ is finite.
\end{thm}

\begin{prf*}
  The inequality $\Gfcgldim \leq \Gwgldim$ holds by
  \fctref{gwgldim}. To prove the opposite inequality one can assume
  that $\Gfcgldim$ is finite, and by \thmref{gfcgf-gldim} it suffices
  to show that $A$ is right weak coherent. Fix an integer $n$ with
  $\Cot\operatorname{-id}_{A}A \leq n$ and $\Gfcgldim\leq n$. We first
  argue that the product of any set of flat-cotorsion $A$-modules has
  finite flat dimension.  This follows from \thmref{8}, since
  $\Gfcgldim\leq n$ implies $\Flat_n \cap \Cot=\Inj_n \cap \Cot$, and
  since $\Inj_n \cap \Cot$ is closed under products, so is
  $\Flat_n \cap \Cot$.

  Taking special cotorsion preenvelopes yields per \prpref{Vid} a
  resolution,
  \begin{equation*}
    0 \lra A \lra C_0 \lra \cdots \lra C_{-n}\lra 0 \:,
  \end{equation*}
  where the $A$-modules $C_i$ are flat-cotorsion.  For any set $X$ it
  yields an exact sequence
  \begin{equation*}
    0 \lra A^X \lra C_0^X \lra \cdots \lra C_{-n}^X \lra 0 \:.
  \end{equation*}
  As $C_i^X$ has finite flat dimension for each $i$, so does
  $A^X$. Hence by \thmcite[4.2]{MCI16}, the ring $A$ is right weak
  coherent. For the last assertion, recall from
  \thmcite[19.2.5(1)]{LMaNDn06} that $\splf$ equals the cotorsion
  global dimension of $A$, see also \rmkref{simple}.
\end{prf*}

In comments on an earlier version of this paper, Ioannis Emmanouil
made the following observation.

\begin{rmk}
  If every flat $A$-module has finite Gorenstein projective dimension,
  then $\Cot\operatorname{-id}_{A}A$ is finite and \thmref{58}
  applies. Indeed, as flat modules are closed under coproducts, there
  is an integer $n$ such that $\Gpd F\le n$ holds for every flat
  $A$-module $F$. To see that $\Cot\operatorname{-id}_{A}A\le n$
  holds, it suffices by \prpref{Vid} to show that
  $\Ext_{A}^{n+1}(F,A)=0$ holds for every flat $A$-module $F$. Let
  $P\qra F$ be a projective resolution. Dimension shifting yields
  $\Ext_A^{n+1}(F,A)\cong \Ext_A^1(\Co[n]{P},A)$, which vanishes
  because $\Co[n]{P}$ is Gorenstein projective.
\end{rmk}

\begin{cor}
  \label{cor:osofsky0}
  If $\splf$ and $\splf[\Aop]$ are both finite, then one has:
  \begin{equation*}
    \Gfcgldim = \Gfcgldim[\Aop] \:.
  \end{equation*}
\end{cor}

\begin{prf*}
  The assertion follows in view of \corcite[2.5]{CET-21a} from
  \thmref{58}.
\end{prf*}

Recall that $A$ is called left $\aleph_0$-noetherian if every left
ideal is countably generated.

\begin{thm}
  \label{thm:a0}
  If $A$ is left $\aleph_0$-noetherian, then the following equality
  holds:
  \begin{equation*}
    \Gfcgldim = \Gwgldim \:.
  \end{equation*}
\end{thm}

\begin{prf*}
  If $\splf$ is finite, the equality holds by \thmref{58}. If $\splf$
  is infinite, then so is $\FPD$ by \prpcite[6]{CUJ70}. From
  \prpcite[2.8]{IEmOTl11} it follows that $\FFD$ is infinite, whence
  $\Gfcgldim$ is infinite by \thmref{8}, and by \fctref{gwgldim} one
  has $\Gfcgldim \le \Gwgldim$ .
\end{prf*}

\begin{cor}
  \label{cor:a0-1}
  If $A$ is $\aleph_0$-noetherian, then the following equality holds:
  \begin{equation*}
    \Gfcgldim = \Gfcgldim[\Aop] \:.
  \end{equation*}
\end{cor}

\begin{prf*}
  The assertion follows in view of \corcite[2.5]{CET-21a} from
  \thmref{a0}.
\end{prf*}

In what remains of this section we derive some further consequences of
the theorems above in combination with \thmref[Theorems~]{Ggldim-new} and
\thmref[]{8}.

  The next
result is also obtained by Wang, Yang, Shao, and Zhang
\thmcite[3.7]{WYSZ}, and it can easily be deduced from a result of El
Maaouy \thmcite[4.16]{REM}. Both arguments are based on the
PGF-modules introduced in \cite{JSrJSt20}; our proof is different. For
commutative rings the inequality is already known to hold from
\corcite[1.2]{DBnNMh10}, see also \rmkcite[5.4(ii)]{IEm12}, and for
right coherent rings it holds by \corcite[3.5]{CET-21a}.

\begin{thm}
  \label{thm:GwGgl_ineq}
  The following inequality holds:
  \begin{equation*}
    \Gwgldim \leq \Ggldim \:.
  \end{equation*}
\end{thm}

\begin{prf*}
  One can assume that $\Ggldim<\infty$ holds, in which case $A$ is
  right weak coherent by \cite[Exa. 4.4(3)]{MCI16}. Now combine
  \corref{1.3} and \thmref{gfcgf-gldim}.
\end{prf*}

The next corollary applies, in particular, if $A$ is left perfect.

\begin{cor}
  \label{cor:perfect}
  If $\FFD=\FPD$ holds, then the following equalities hold:
  \begin{equation*}
    \Gfcgldim = \Gwgldim = \Ggldim \:.
  \end{equation*}
\end{cor}

\begin{prf*}
  If $\FFD=\FPD$ holds, then \thmref[Theorems~]{Ggldim-new} and
  \thmref[]{8} yield
  \begin{equation*}
    \Ggldim=\max\{\sfli,\FPD\}=\max\{\sfli,\FFD\}=\Gfcgldim\:.
  \end{equation*}
  If this quantity is infinite, then \fctref{gwgldim} implies that
  $\Gwgldim$ is infinite. If the quantity is finite, then
  \thmref{GwGgl_ineq} implies that $\Gwgldim$ is finite and thus
  agrees with $\Gfcgldim$ per \fctref{gwgldim}.
\end{prf*}

If $A$ is noetherian it is known from \corcite[6.11]{ABl00} that
$\Ggldim$ and $\Ggldim[\Aop]$ are simultaneously finite. In the
terminology from \dfncite[6.8]{ABl00} or \dfncite[VII.2.5]{ABlIRt07}
the final two theorems say that the Gorenstein property is symmetric
for $\aleph_0$-noetherian rings and rings of cardinality $\aleph_n$
for some $n \ge 0$.

The next inequality is a counterpart in Gorenstein homological algebra
to \thmcite[1]{CUJ66}, and the subsequent corollary is a counterpart
to \corcite[1]{CUJ66}.

\begin{thm}
  \label{thm:GwGl-1}
  If $A$ is left $\aleph_0$-noetherian, then the next inequality
  holds:
  \begin{equation*}
    \Ggldim \le \Gwgldim + 1 \:.
  \end{equation*}
\end{thm}

\begin{prf*}
  One can assume that $\Gwgldim=n<\infty$ holds. Thus \fctref{gwgldim}
  yields $\Gfcgldim=n$. From \thmref{8} it now follows that the equalities
  $\sfli=\FFD= n$ hold. Since \prpcite[2.8]{IEmOTl11} yields
  $\FPD \leq n+1$, the desired
  inequality $\Ggldim\leq n+1$ follows from \thmref{Ggldim-new}.
\end{prf*}

\begin{cor}
  \label{cor:GwGl-2}
  If $A$ is $\aleph_0$-noetherian, then $\Ggldim$ is finite if and
  only if $\Ggldim[\Aop]$ is finite, and in that case the following
  inequality holds:
  \begin{equation*}
    |\Ggldim - \Ggldim[\Aop]| \le 1 \:.
  \end{equation*}
\end{cor}

\begin{prf*}
  This follows from \thmref{GwGl-1} via \thmref{GwGgl_ineq} and
  \corcite[2.5]{CET-21a}.
\end{prf*}

\begin{rmk}
  Under the a priori stronger assumption that both $\Ggldim$ and
  $\Ggldim[\Aop]$ are finite, the inequality in the previous corollary
  was noticed in \rmkcite[5.4(iii)]{IEm12}. One could alternatively
  obtain \corref{GwGl-2} by combining \thmref{GwGgl_ineq} with
  \thmcite[2.9]{IEmOTl11} and a result of Dalezios and Emmanouil
  \prpcite[5.3]{GDlIEm}.
\end{rmk}

The last assertion in the next result is a counterpart in Gorenstein
homological algebra to a result of Osofsky \corcite[1.6]{BLO68a}. As a
countable ring is $\aleph_0$-noetherian, the case $n=0$ is a special
instance of \corref{GwGl-2}.

\begin{thm}
  \label{thm:Osofsky}
  If $\splf$ and $\splf[\Aop]$ are both finite, then $\Ggldim$ is finite if and
  only if $\Ggldim[\Aop]$ is finite, and in that case  there is an inequality:
  \begin{equation*}
    |\Ggldim - \Ggldim[\Aop]| \le \max\{\splf,\splf[\Aop] \} \:.
  \end{equation*}
  In particular, if $A$ has cardinality $\aleph_n$ for an integer
  $n \ge 0$, then one has:
  \begin{equation*}
    |\Ggldim - \Ggldim[\Aop]| \le n + 1 \:.
  \end{equation*}
\end{thm}

\begin{prf*}
  By \corref{1.3} and \corref{osofsky0} there are inequalities
  \begin{align*}
    \Ggldim & \le \Gfcgldim + \splf \\
            & = \Gfcgldim[\Aop] + \splf
              \le \Ggldim[\Aop] + \splf \:,
  \end{align*}
  and by symmetry one has $\Ggldim[\Aop] \le \Ggldim +
  \splf[\Aop]$. This proves the first inequality. For the last
  assertion recall that by a result of Simson \cite{DSm74} a ring of
  cardinality $\aleph_n$ has $\splf \leq n+1$ and
  $\splf[\Aop]\leq n+1$.
\end{prf*}

\section*{Acknowledgment}

\noindent 
We thank Ioannis Emmanouil for valuable comments on a draft of
this~paper. Thanks are also due to the anonymous refereee for
suggestions that improved the exposition.



\def\soft#1{\leavevmode\setbox0=\hbox{h}\dimen7=\ht0\advance \dimen7
  by-1ex\relax\if t#1\relax\rlap{\raise.6\dimen7
  \hbox{\kern.3ex\char'47}}#1\relax\else\if T#1\relax
  \rlap{\raise.5\dimen7\hbox{\kern1.3ex\char'47}}#1\relax \else\if
  d#1\relax\rlap{\raise.5\dimen7\hbox{\kern.9ex \char'47}}#1\relax\else\if
  D#1\relax\rlap{\raise.5\dimen7 \hbox{\kern1.4ex\char'47}}#1\relax\else\if
  l#1\relax \rlap{\raise.5\dimen7\hbox{\kern.4ex\char'47}}#1\relax \else\if
  L#1\relax\rlap{\raise.5\dimen7\hbox{\kern.7ex
  \char'47}}#1\relax\else\message{accent \string\soft \space #1 not
  defined!}#1\relax\fi\fi\fi\fi\fi\fi}
  \providecommand{\MR}[1]{\mbox{\href{http://www.ams.org/mathscinet-getitem?mr=#1}{#1}}}
  \renewcommand{\MR}[1]{\mbox{\href{http://www.ams.org/mathscinet-getitem?mr=#1}{#1}}}
  \providecommand{\arxiv}[2][AC]{\mbox{\href{http://arxiv.org/abs/#2}{\sf
  arXiv:#2 [math.#1]}}} \def\cprime{$'$}
\providecommand{\bysame}{\leavevmode\hbox to3em{\hrulefill}\thinspace}
\providecommand{\MR}{\relax\ifhmode\unskip\space\fi MR }
\providecommand{\MRhref}[2]{%
  \href{http://www.ams.org/mathscinet-getitem?mr=#1}{#2}
}
\providecommand{\href}[2]{#2}

\end{document}